\input amssym.tex
\parindent=0cm
\hoffset=.7truecm
\voffset=-.5truecm
\hsize=15.5truecm
\vsize=24.5truecm
\baselineskip 15pt
\input graphicx 

\def\coordinates{{\includegraphics[width=.3\hsize]{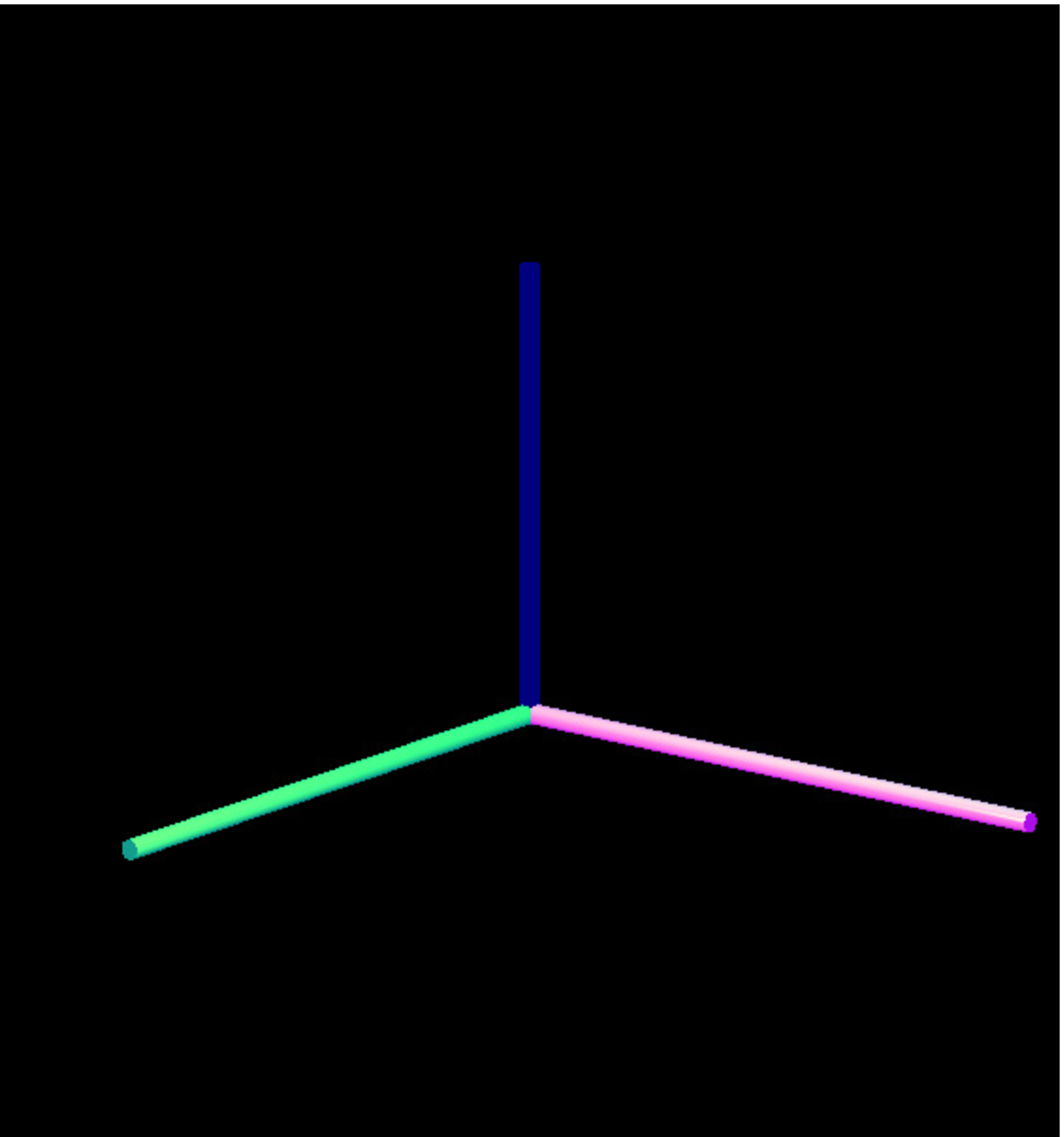}}}
\def\clipA{{\includegraphics[width=.45\hsize]{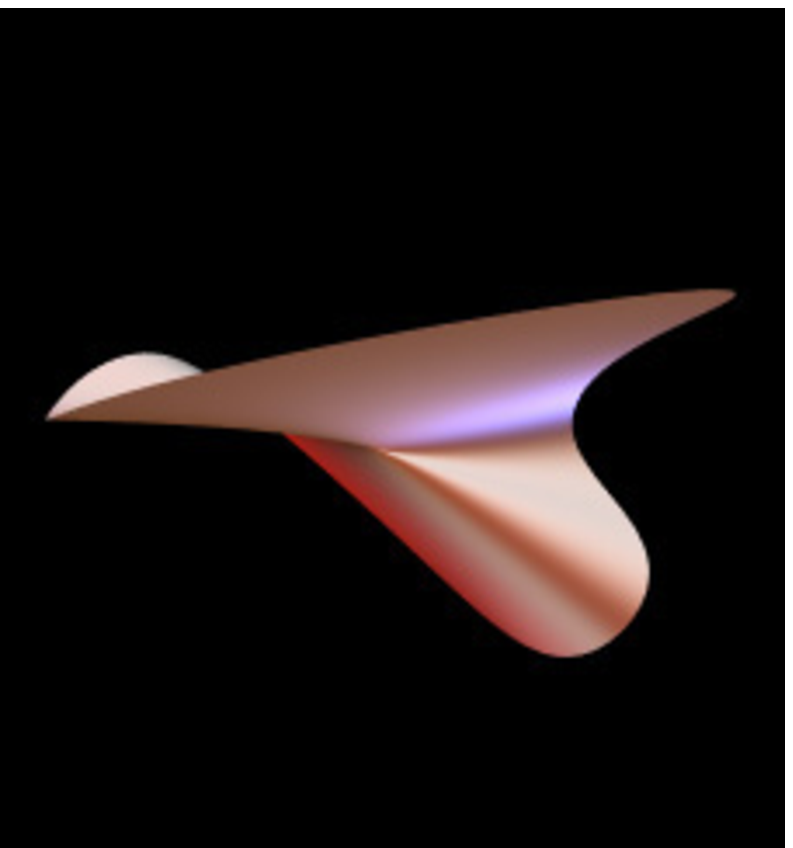}}}
\def\columpiusA{{\includegraphics[width=.45\hsize]{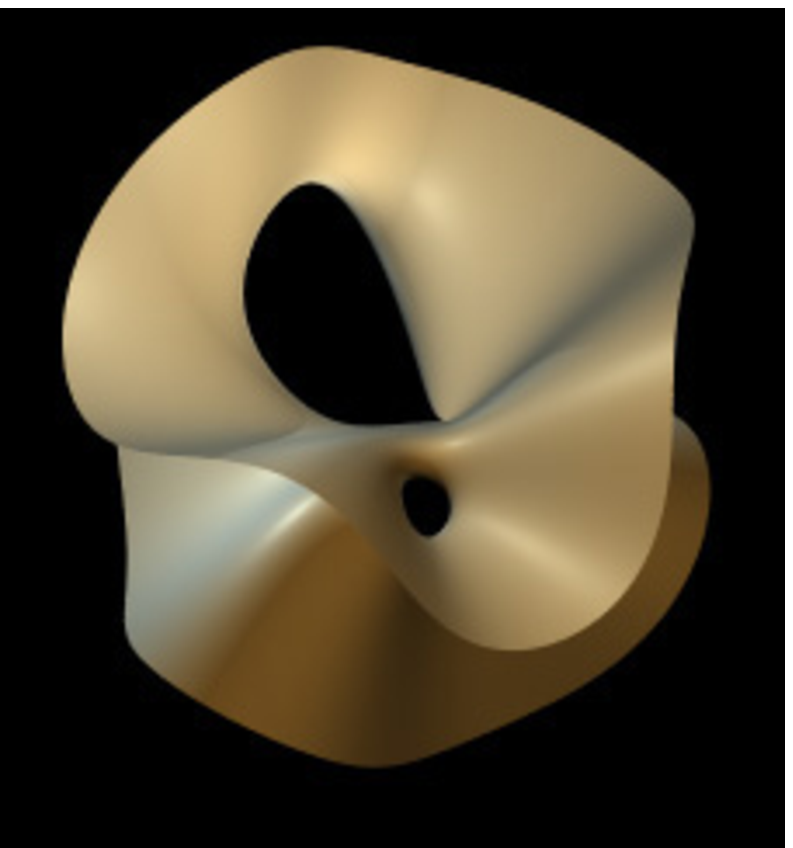}}}
\def\croissantA{{\includegraphics[width=.45\hsize]{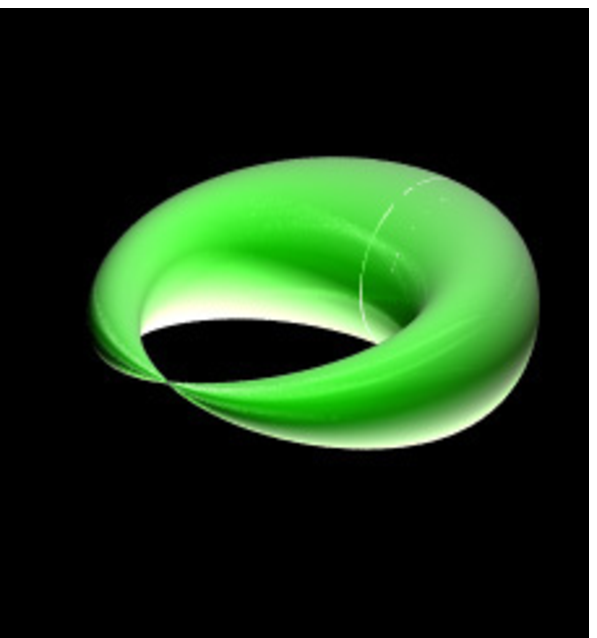}}}
\def\cubeA{{\includegraphics[width=.45\hsize]{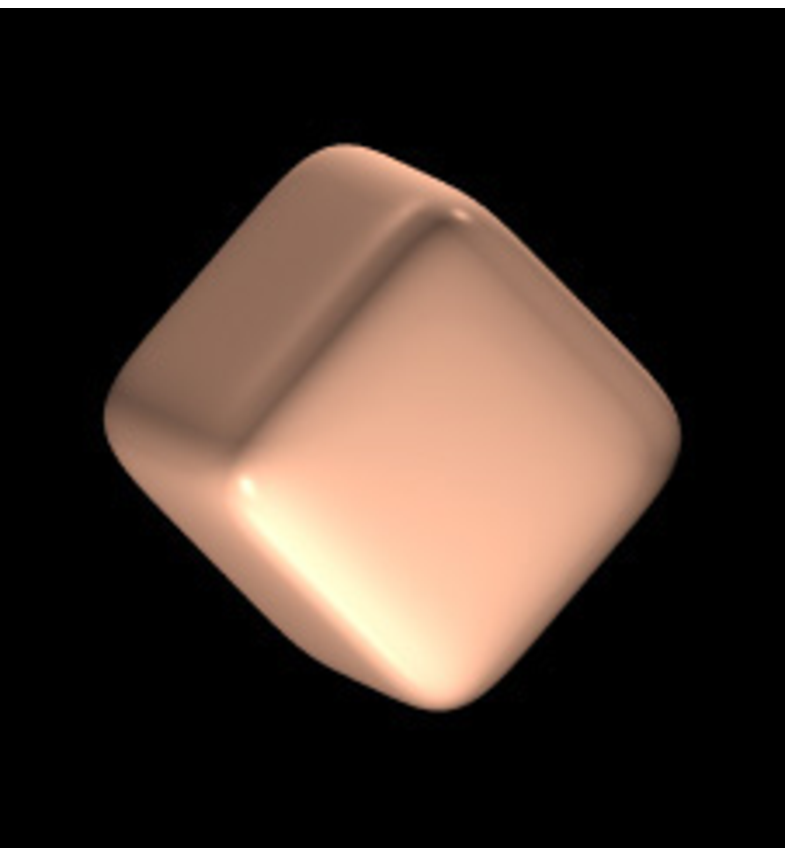}}}
\def\daisyA{{\includegraphics[width=.45\hsize]{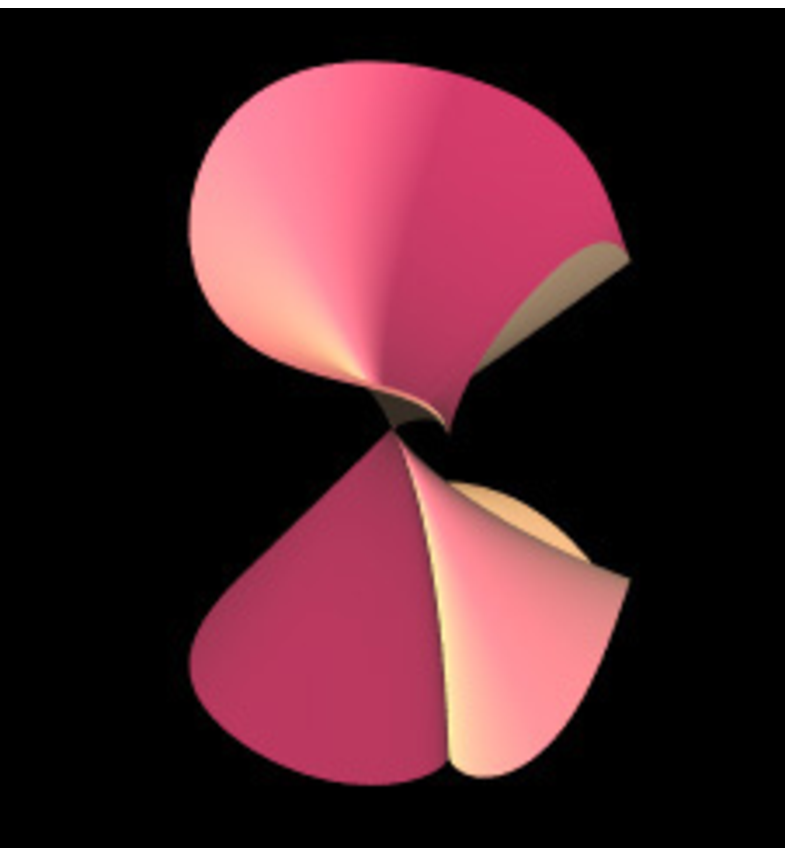}}}
\def\dattelA{{\includegraphics[width=.45\hsize]{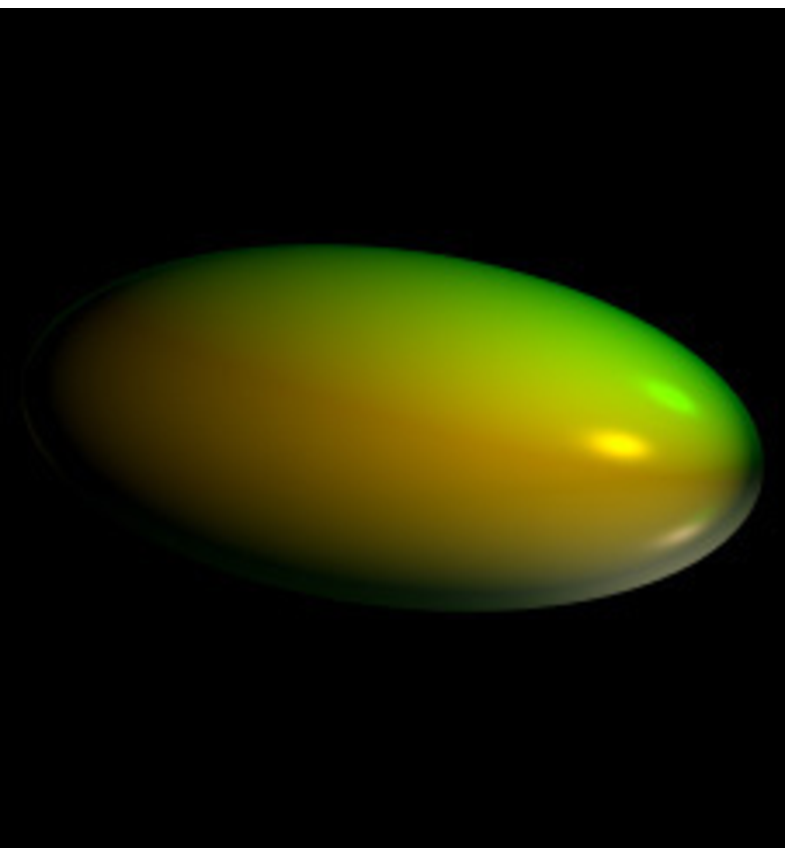}}}
\def\distelA{{\includegraphics[width=.45\hsize]{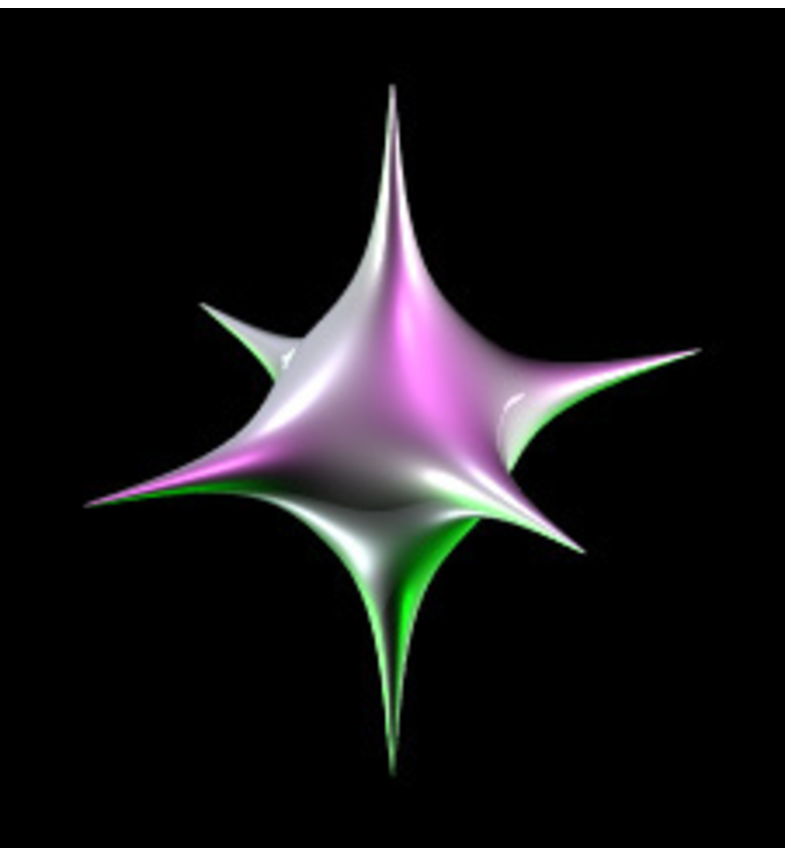}}}
\def\dromoA{{\includegraphics[width=.45\hsize]{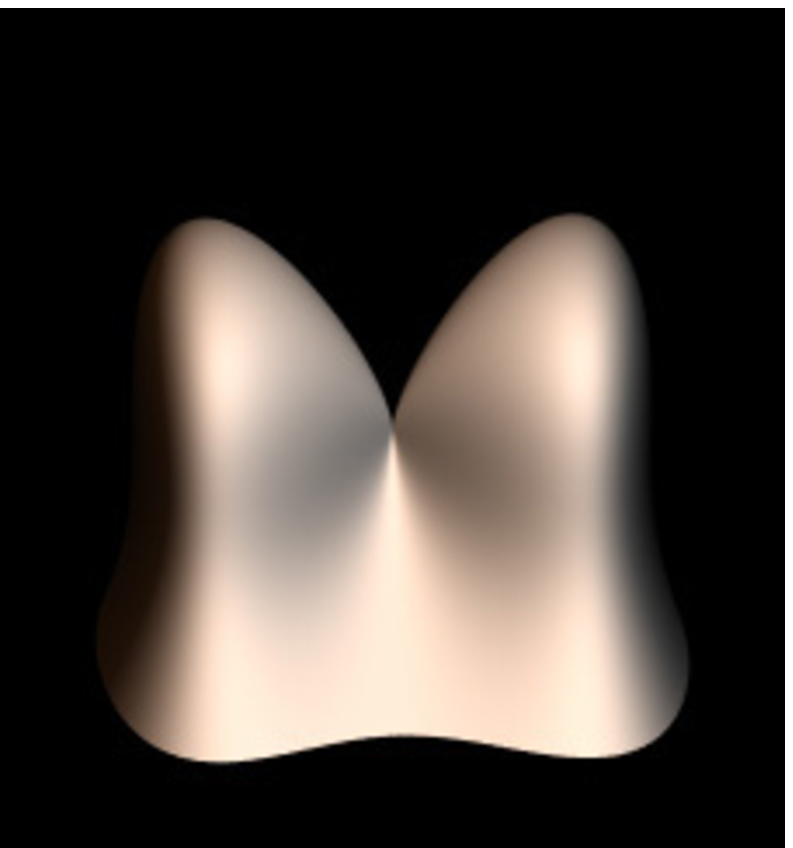}}}
\def\durchblickA{{\includegraphics[width=.45\hsize]{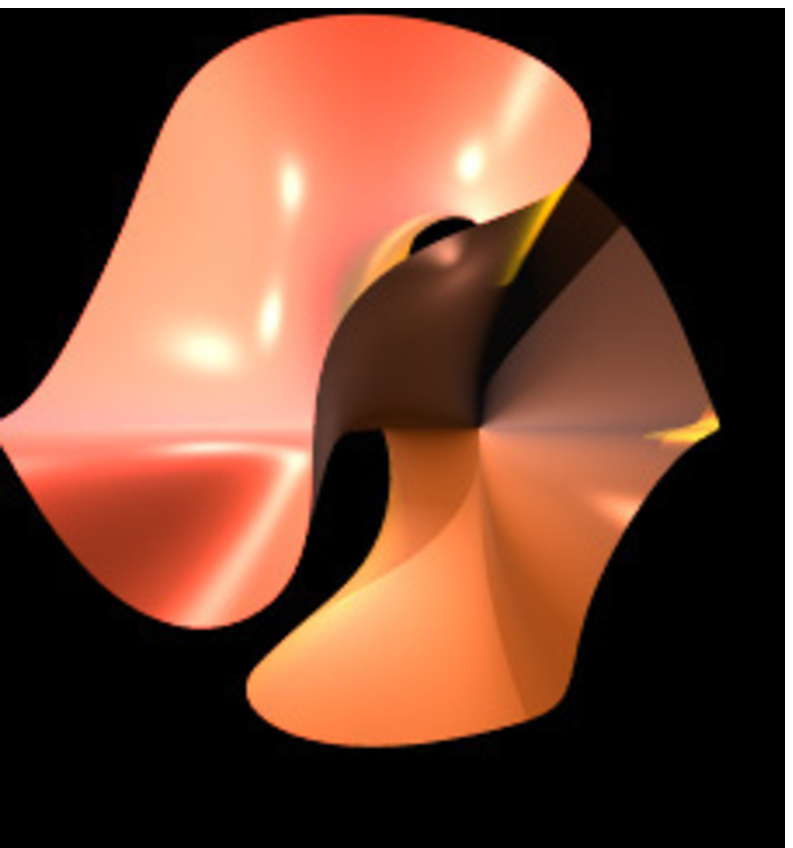}}}
\def\flirtA{{\includegraphics[width=.45\hsize]{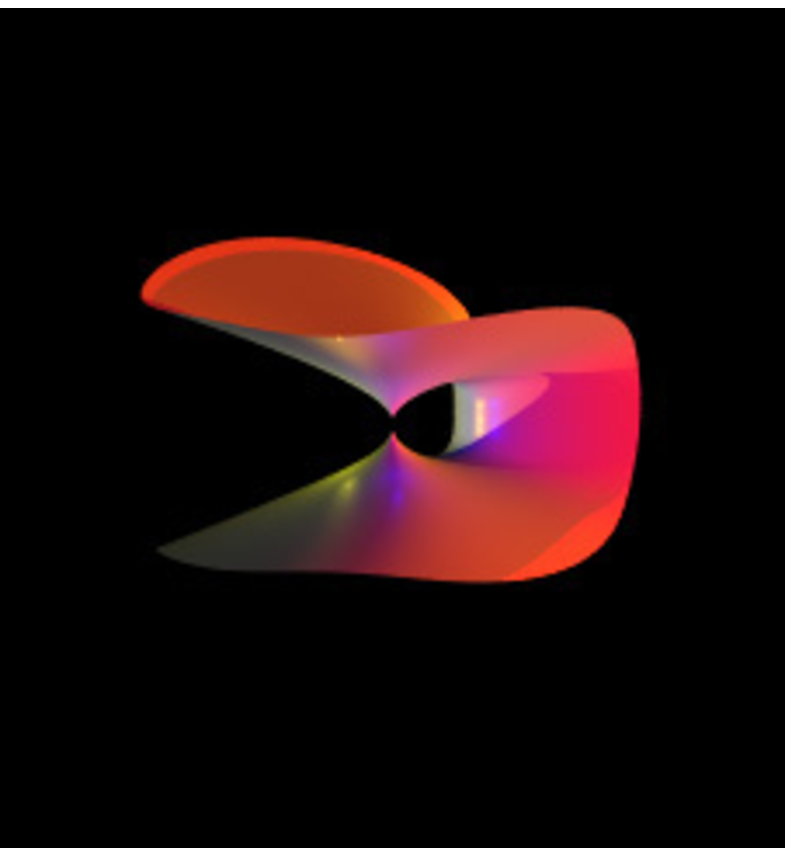}}}
\def\helixA{{\includegraphics[width=.45\hsize]{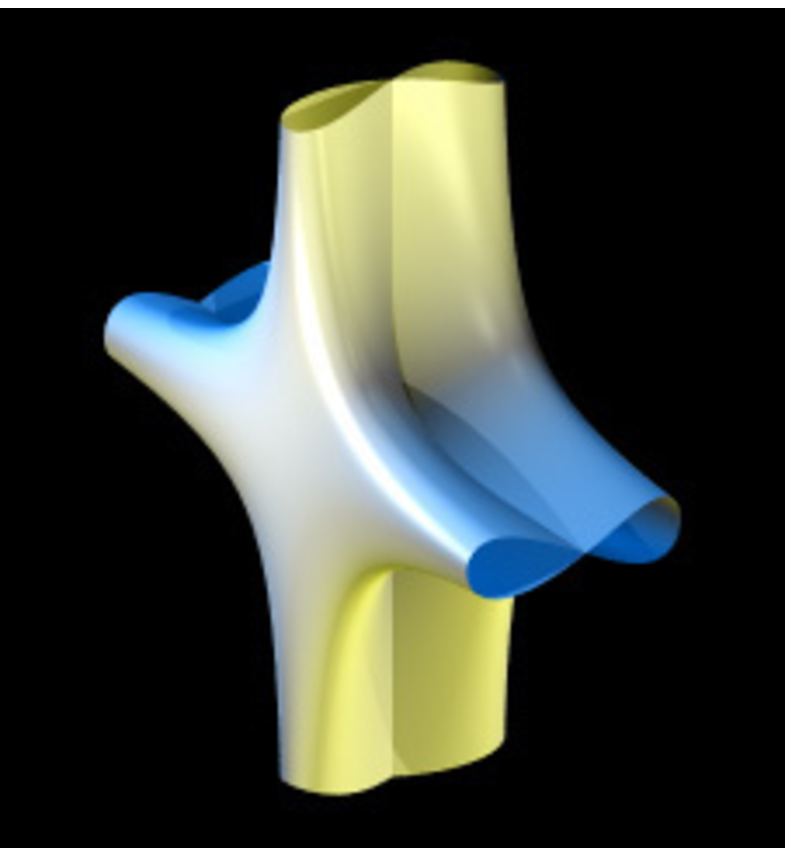}}}
\def\leopoldA{{\includegraphics[width=.45\hsize]{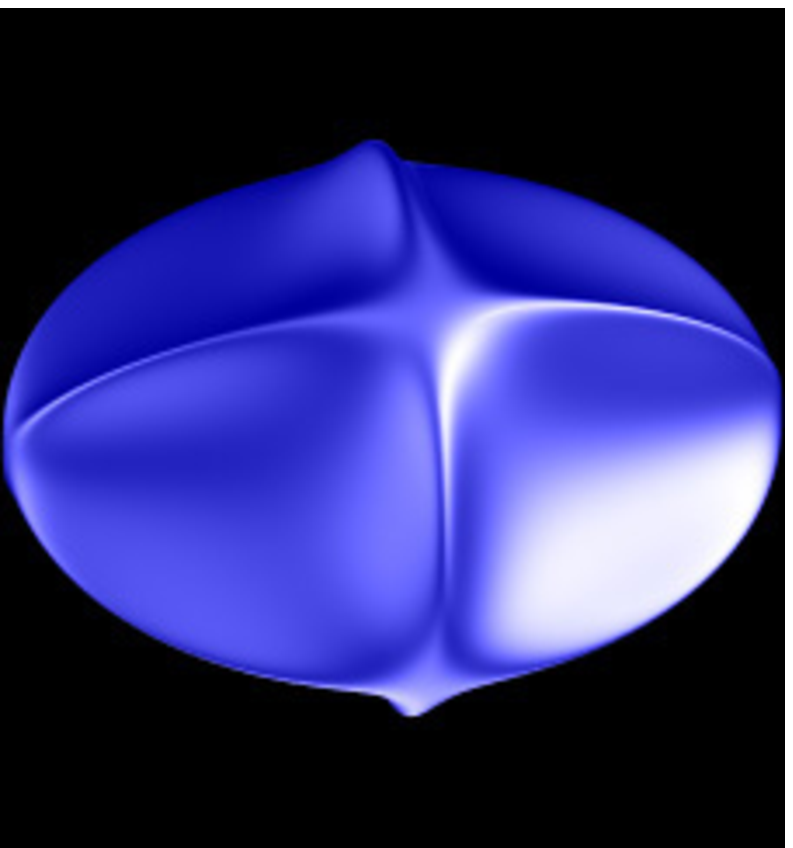}}}
\def\lilieA{{\includegraphics[width=.45\hsize]{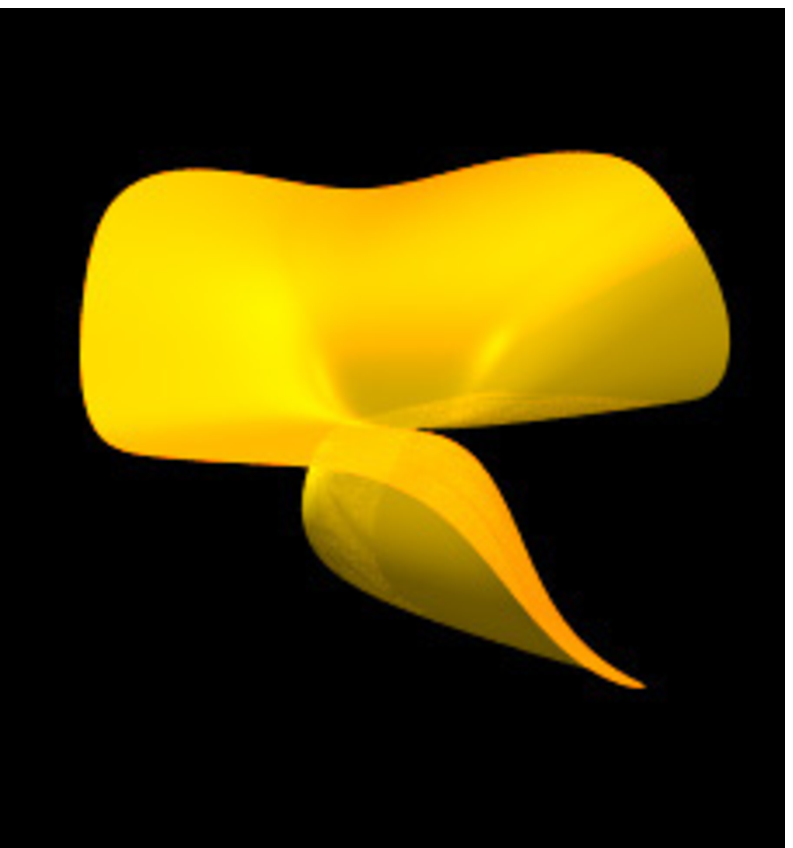}}}
\def\nadeloehrA{{\includegraphics[width=.45\hsize]{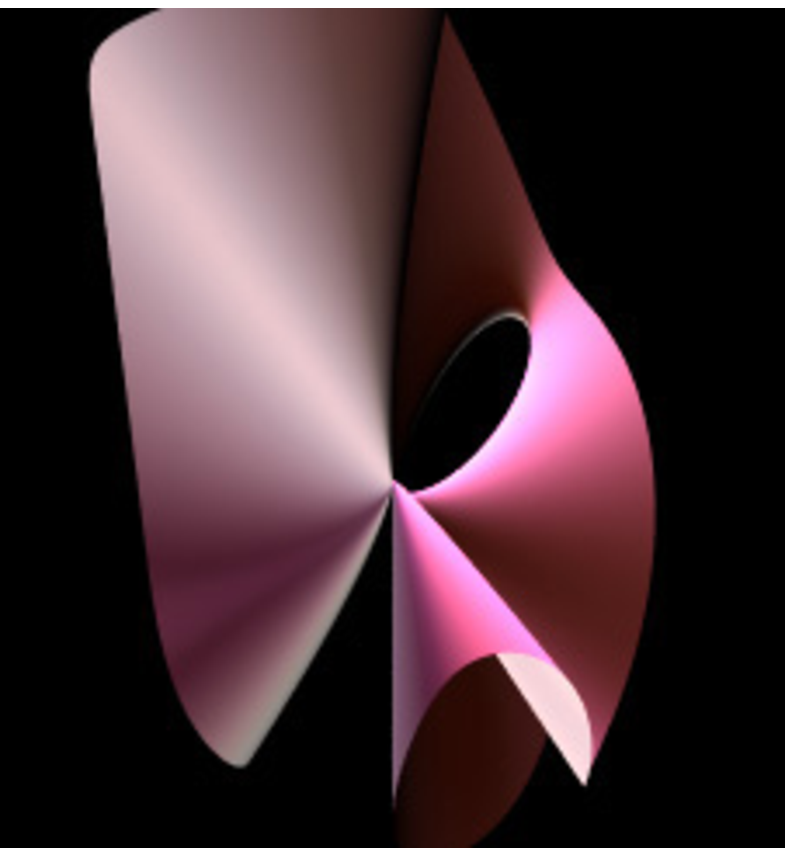}}}
\def\piratA{{\includegraphics[width=.45\hsize]{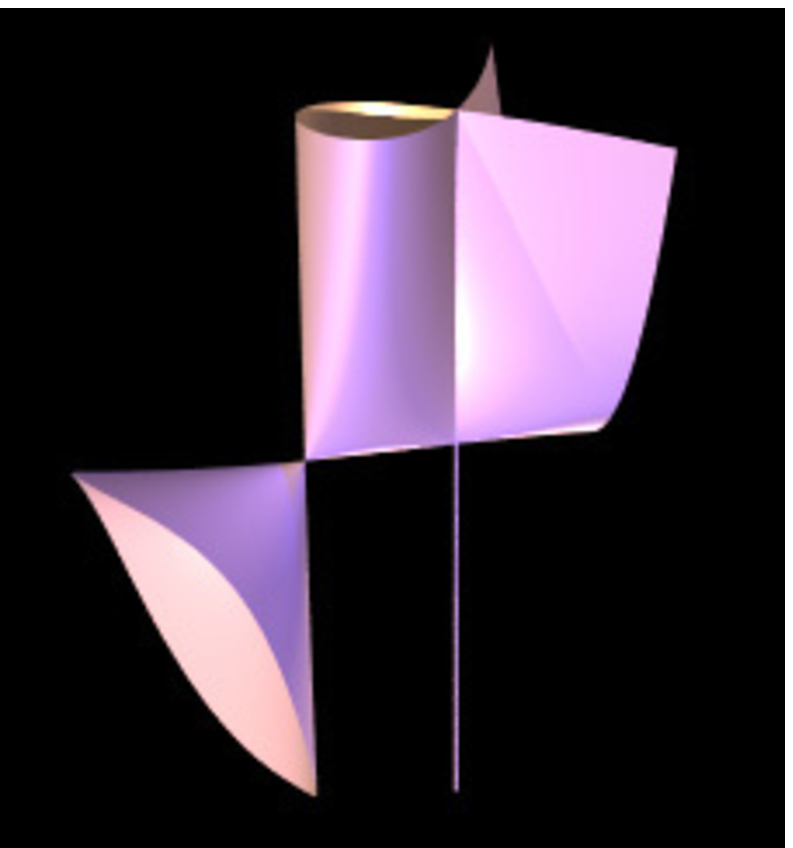}}}
\def\quasteA{{\includegraphics[width=.45\hsize]{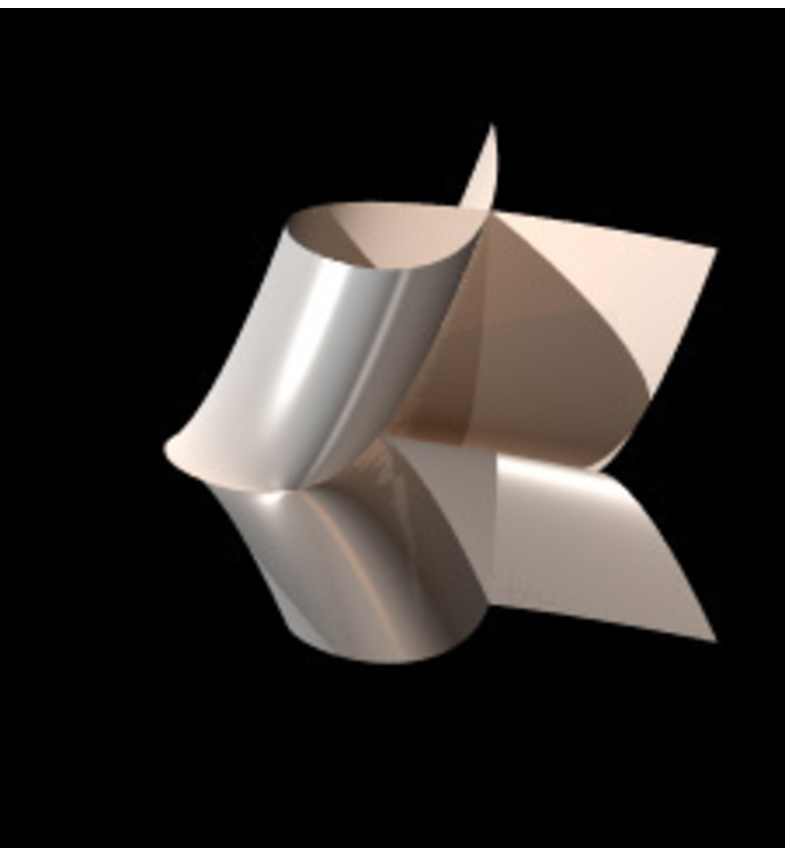}}}
\def\schneeflockeA{{\includegraphics[width=.45\hsize]{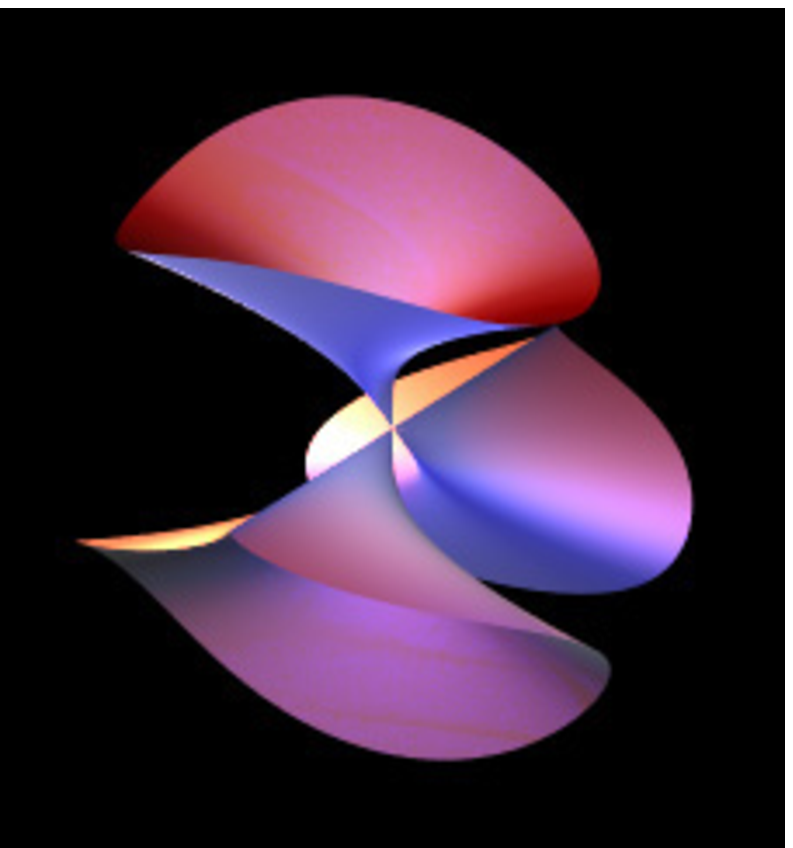}}}
\def\suessA{{\includegraphics[width=.45\hsize]{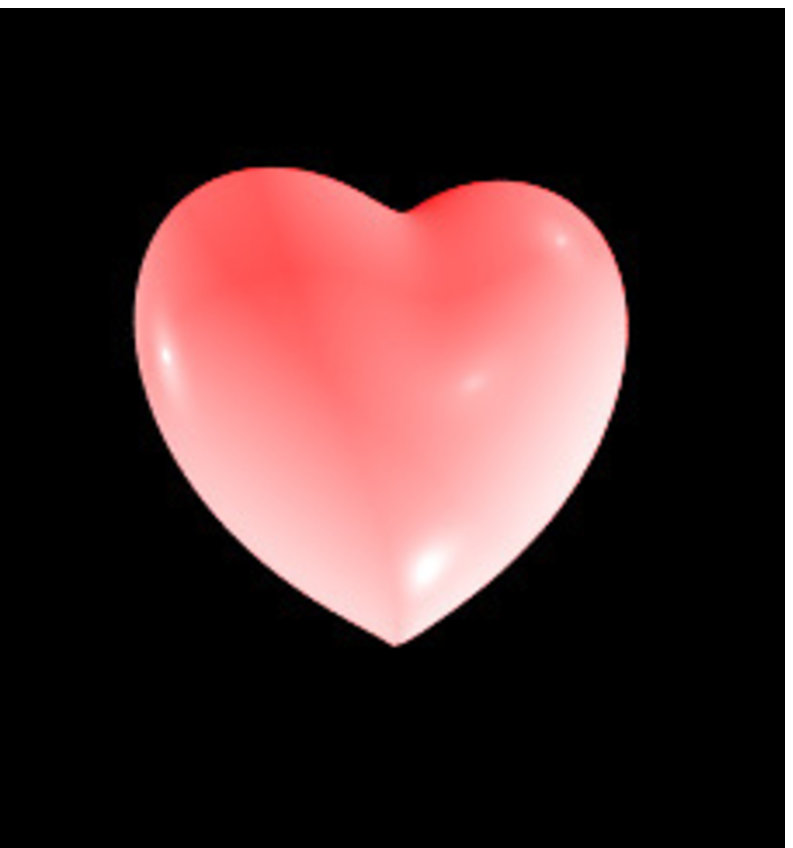}}}
\def\tanzA{{\includegraphics[width=.45\hsize]{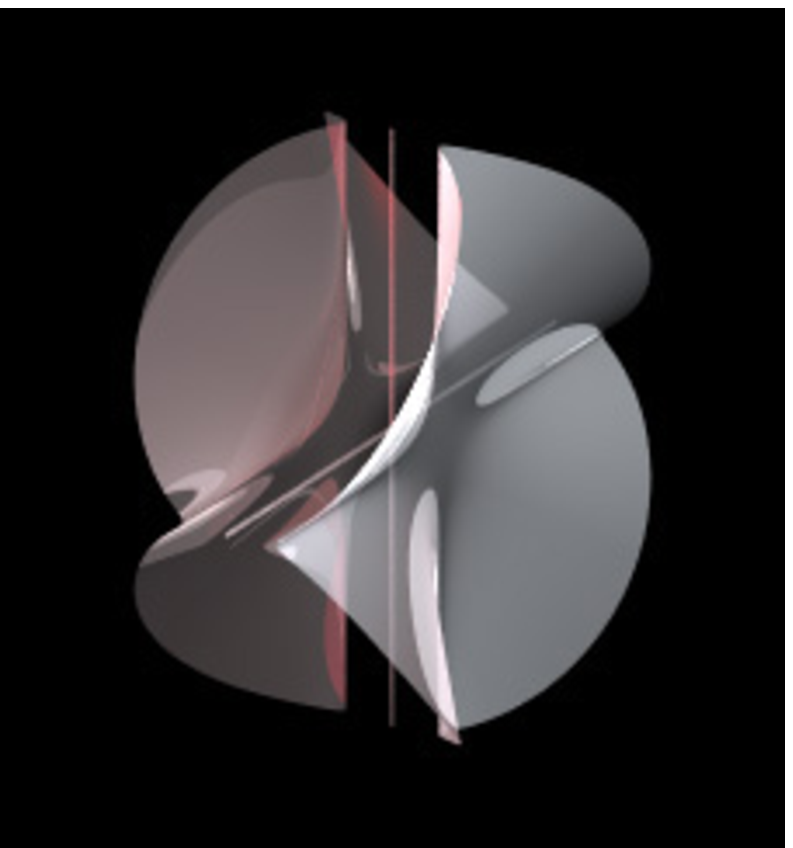}}}
\def\taubeA{{\includegraphics[width=.45\hsize]{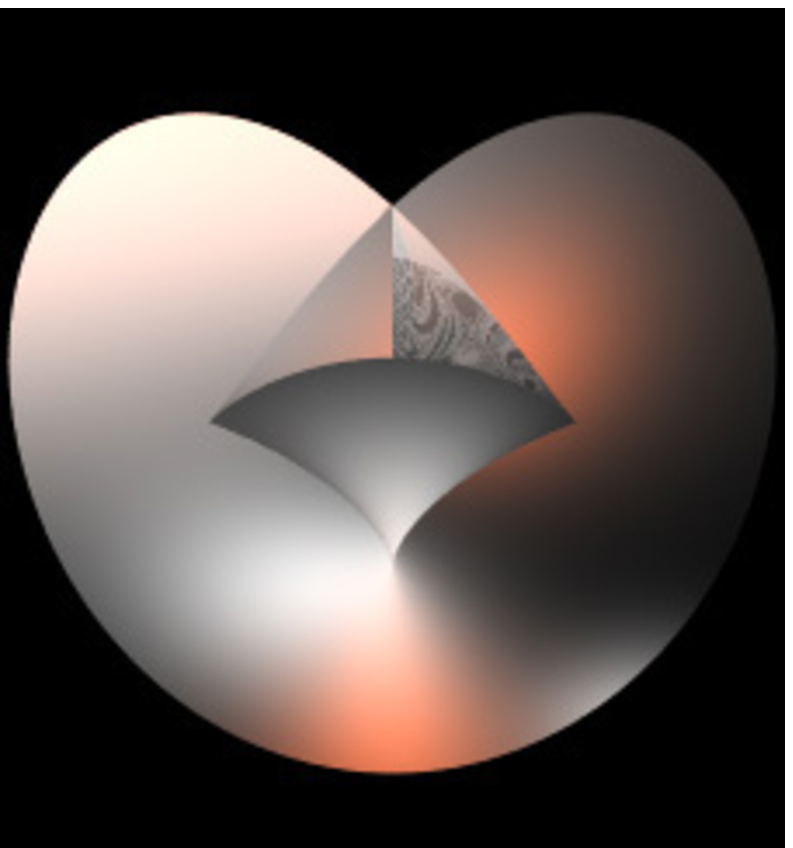}}}
\def\tuelleA{{\includegraphics[width=.45\hsize]{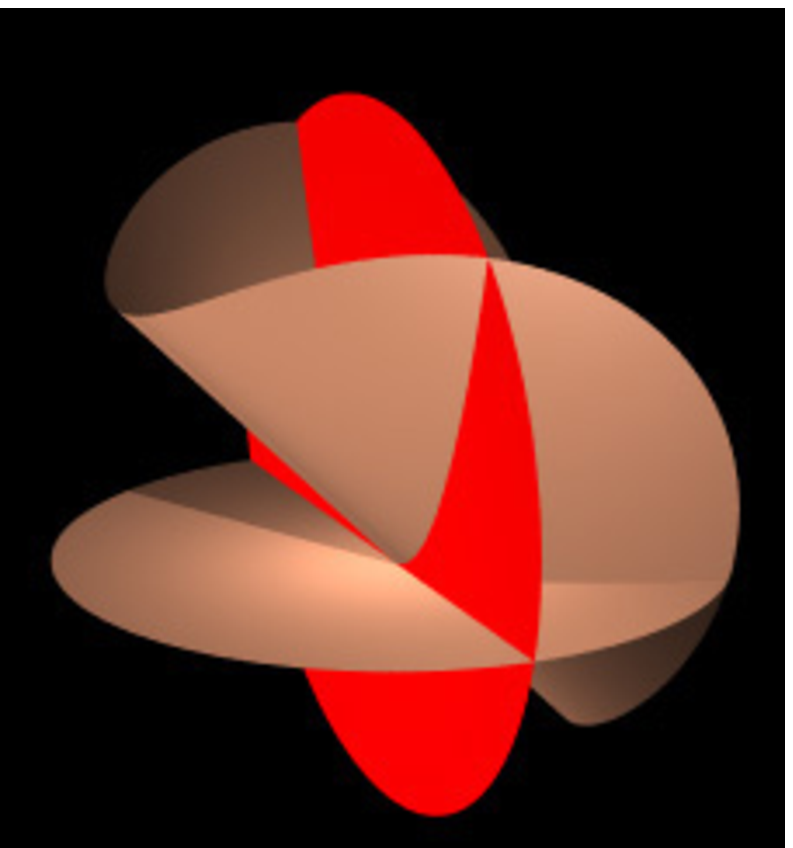}}}

\def\visavisA{{\includegraphics[width=.45\hsize]{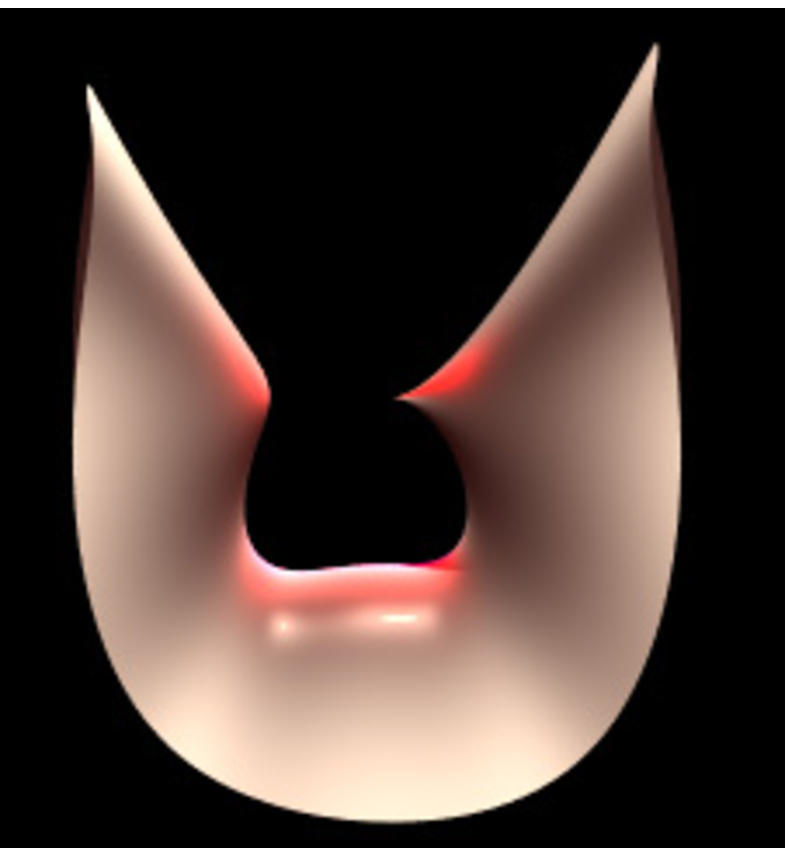}}}
\def\zeppelinA{{\includegraphics[width=.45\hsize]{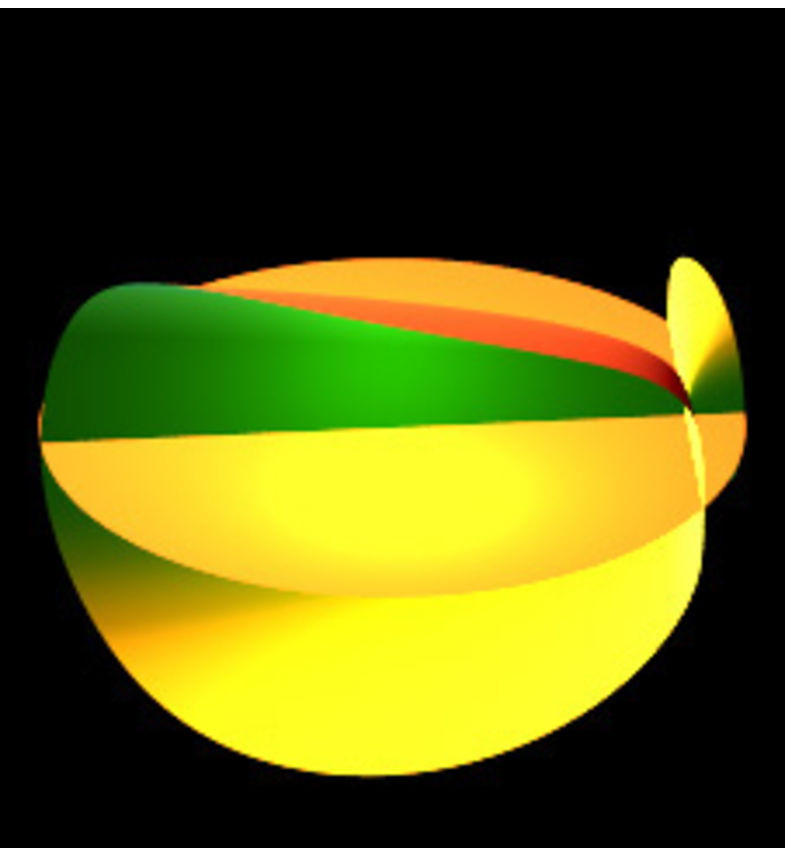}}}
\def\zitrusA{{\includegraphics[width=.45\hsize]{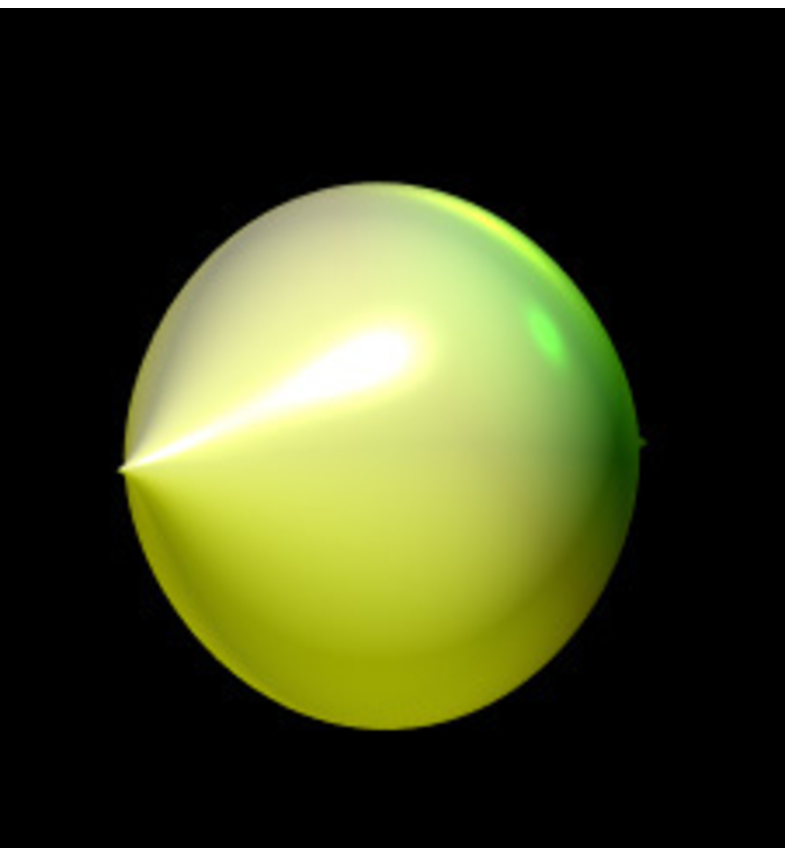}}}

\def\clipB{{\includegraphics[width=.45\hsize]{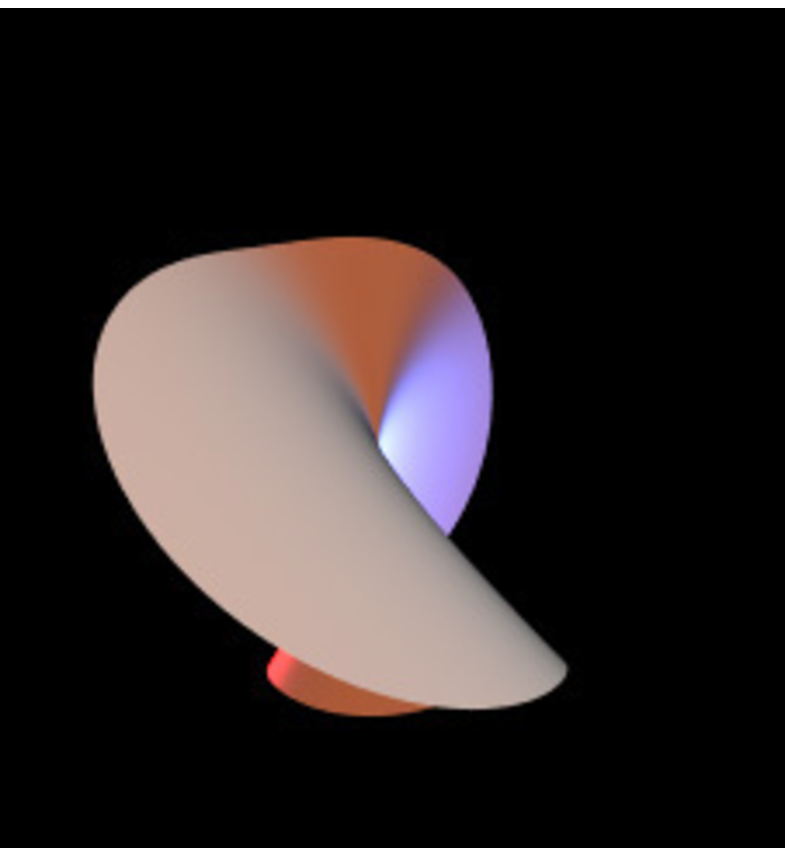}}}
\def\columpiusB{{\includegraphics[width=.45\hsize]{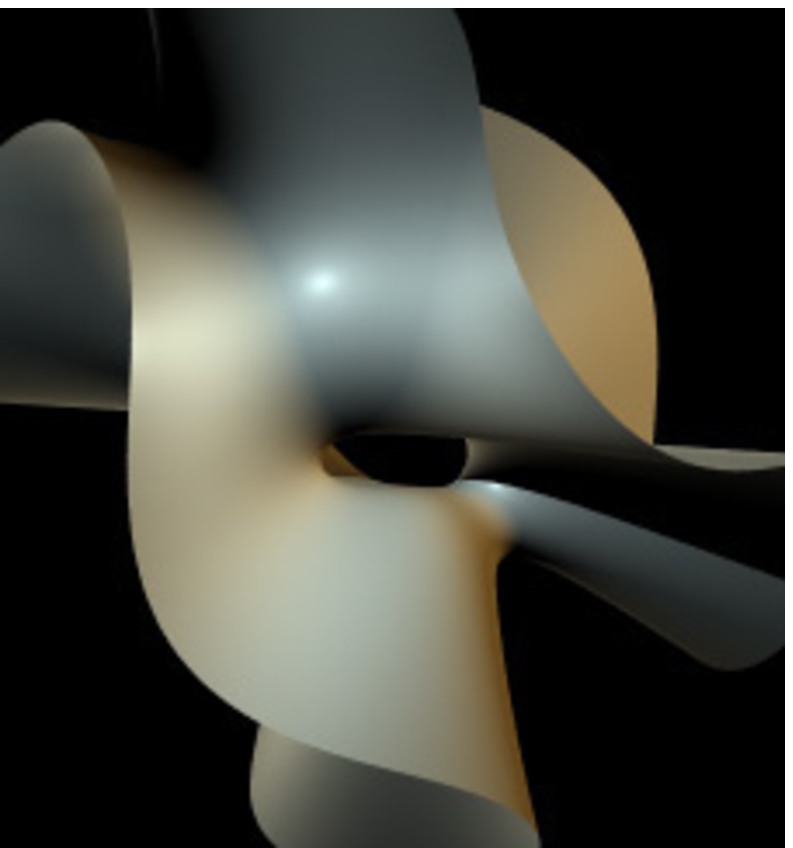}}}
\def\croissantB{{\includegraphics[width=.45\hsize]{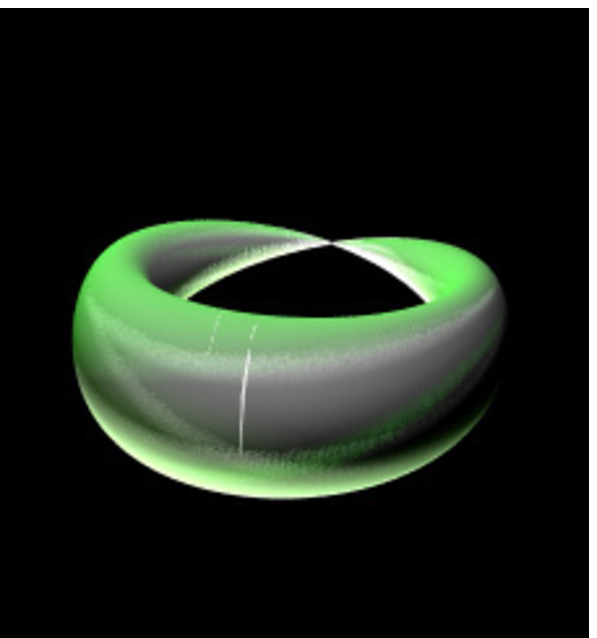}}}
\def\cubeB{{\includegraphics[width=.45\hsize]{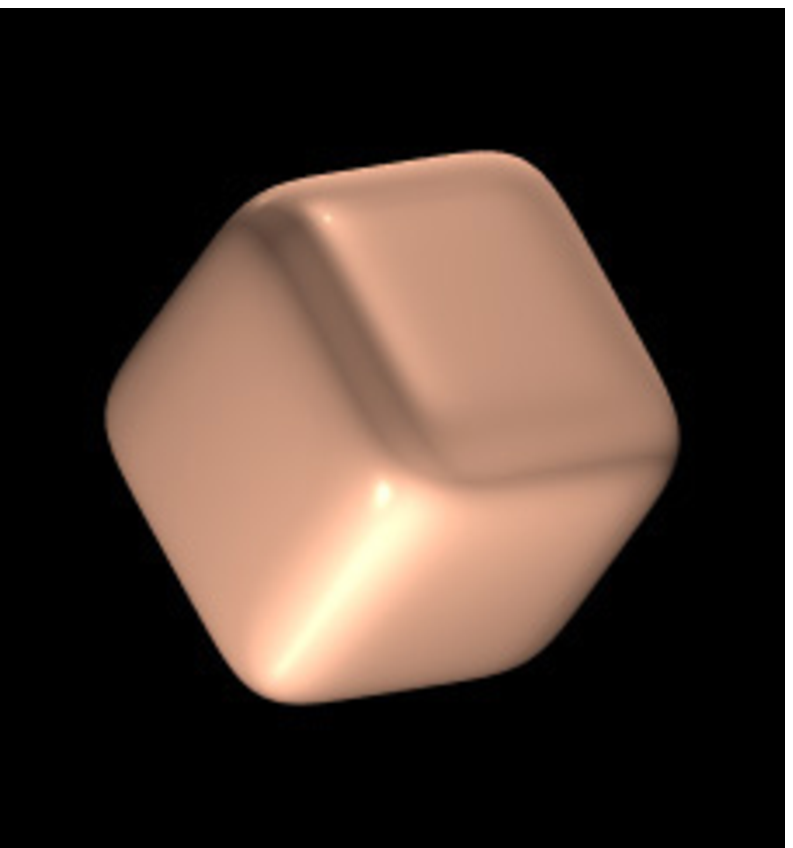}}}
\def\daisyB{{\includegraphics[width=.45\hsize]{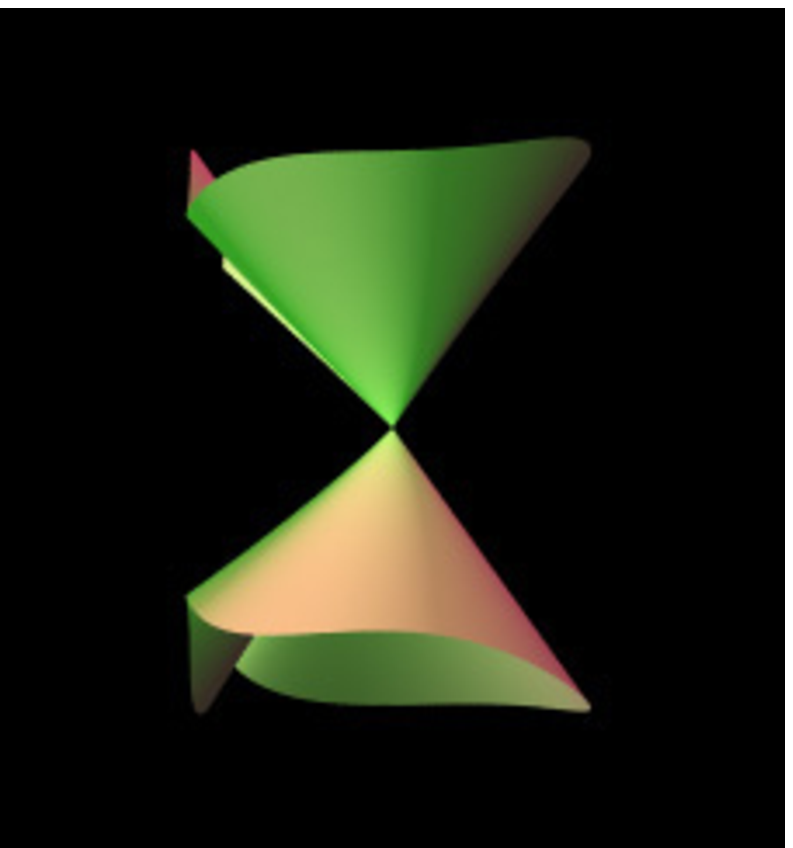}}}
\def\dattelB{{\includegraphics[width=.45\hsize]{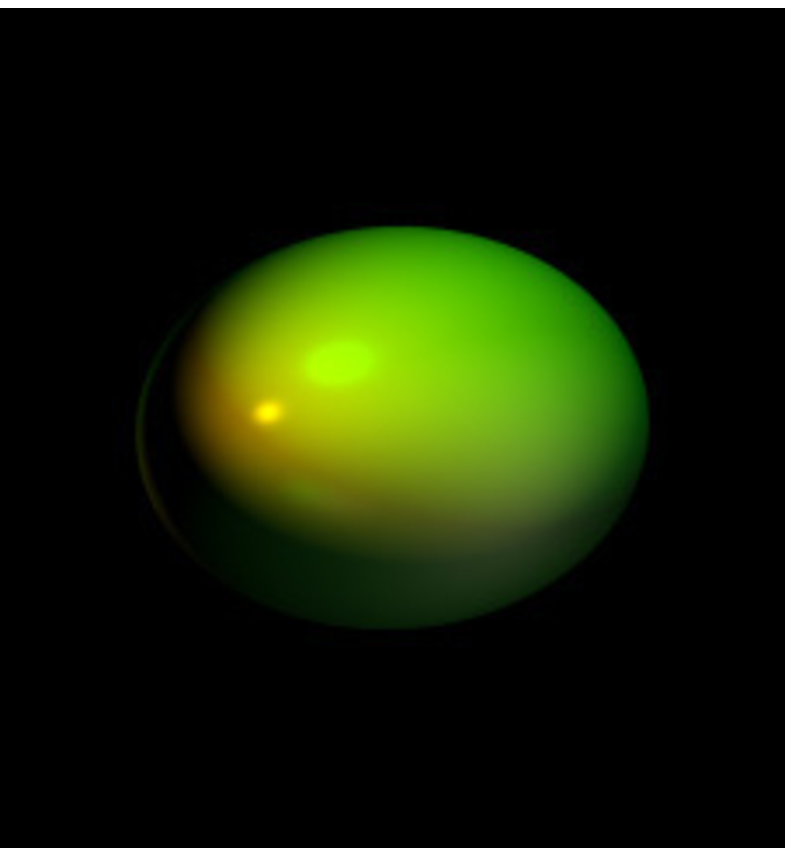}}}
\def\distelB{{\includegraphics[width=.45\hsize]{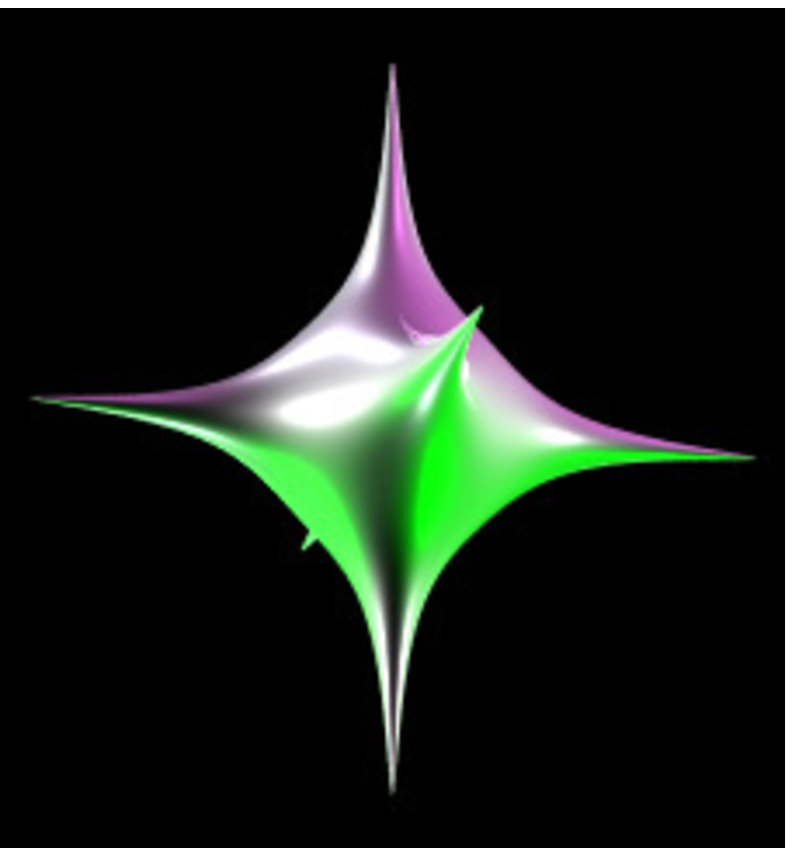}}}
\def\dromoB{{\includegraphics[width=.45\hsize]{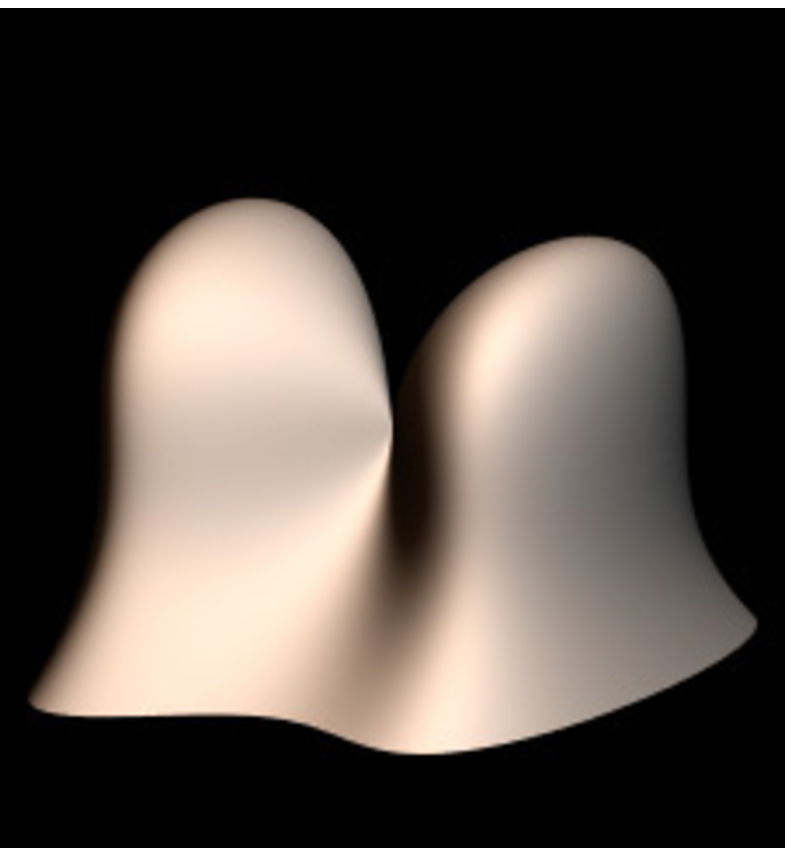}}}
\def\durchblickB{{\includegraphics[width=.45\hsize]{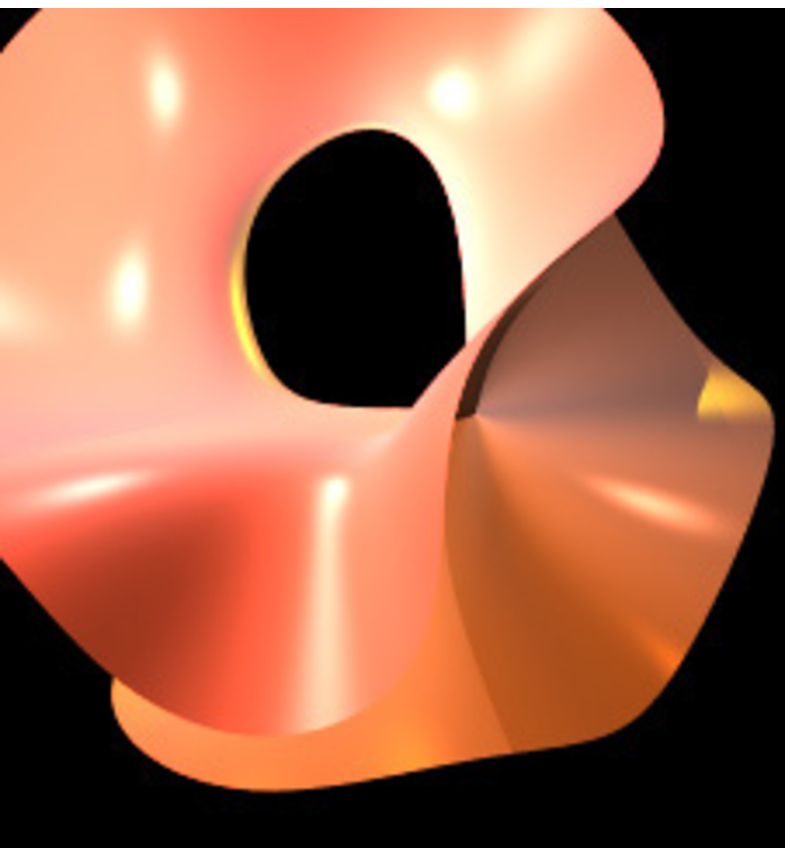}}}
\def\flirtB{{\includegraphics[width=.45\hsize]{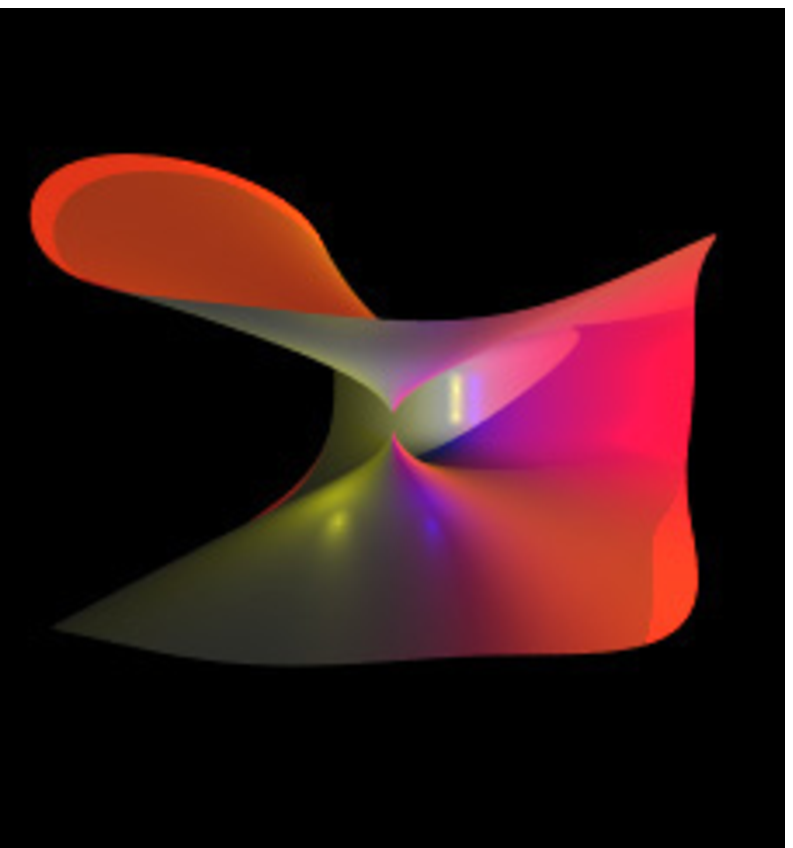}}}
\def\helixB{{\includegraphics[width=.45\hsize]{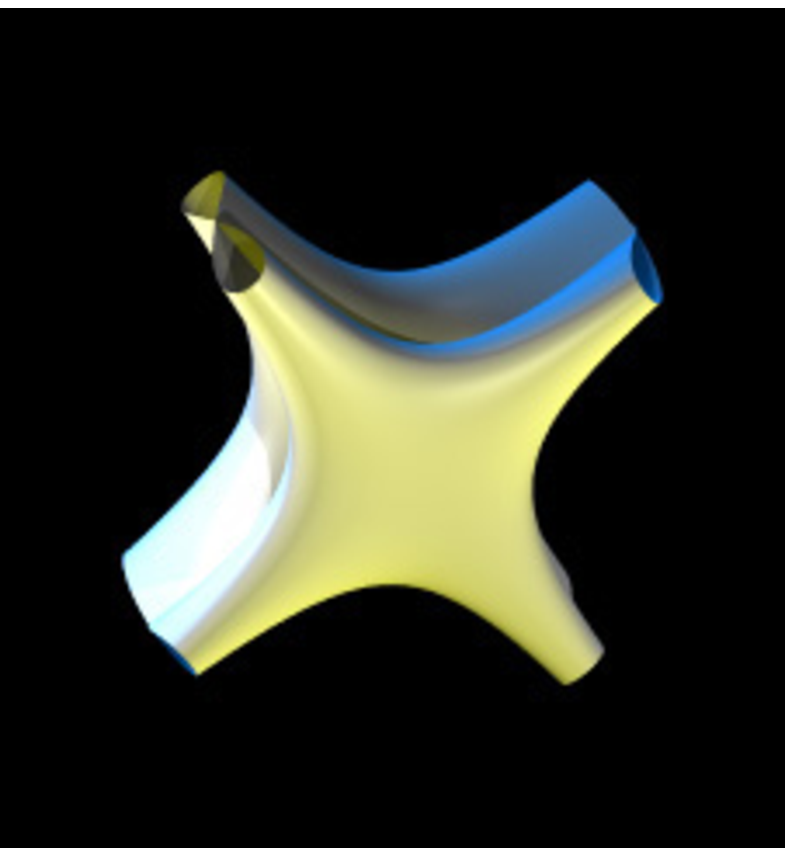}}}
\def\leopoldB{{\includegraphics[width=.45\hsize]{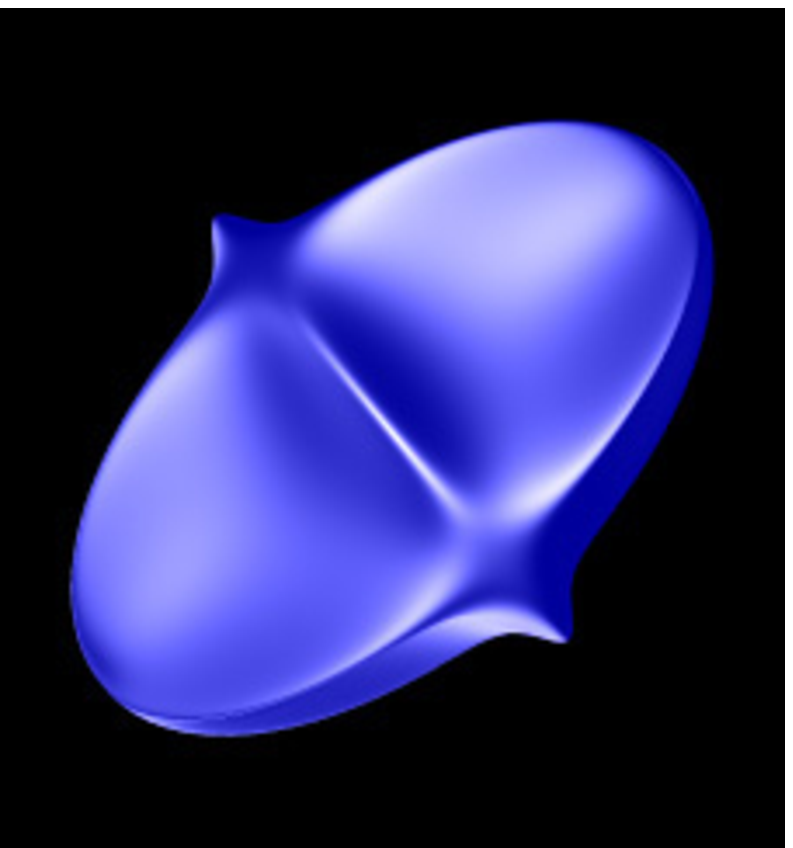}}}
\def\lilieB{{\includegraphics[width=.45\hsize]{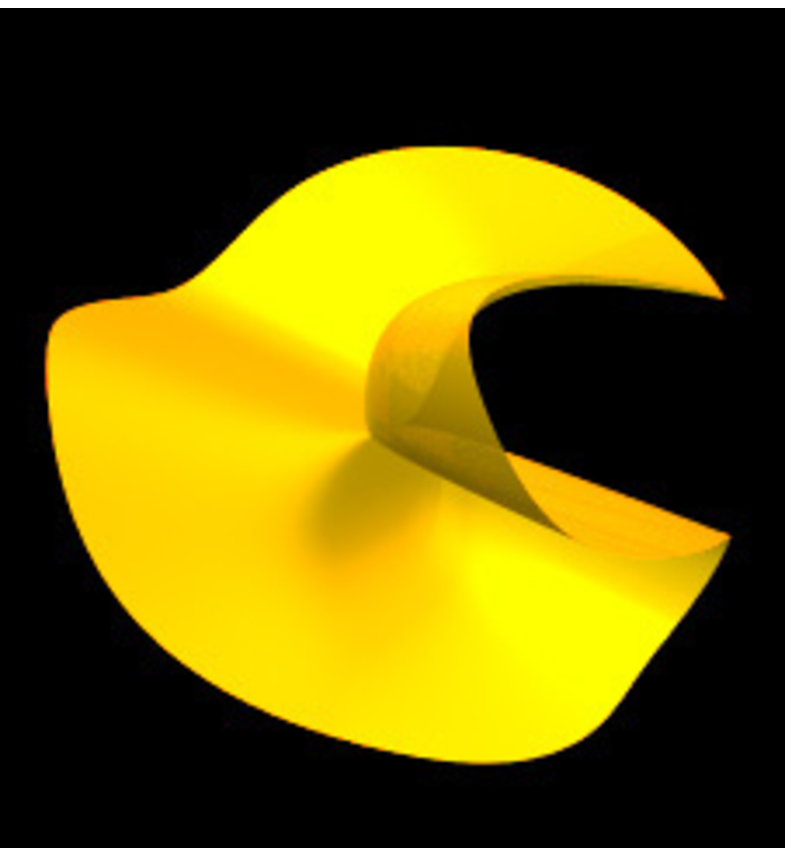}}}
\def\nadeloehrB{{\includegraphics[width=.45\hsize]{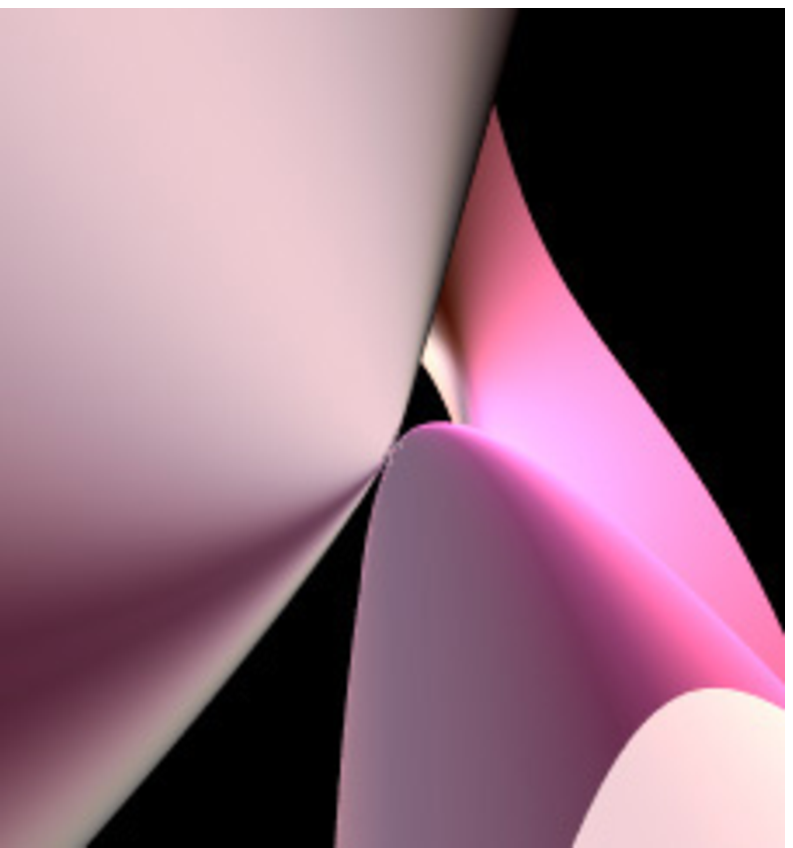}}}
\def\piratB{{\includegraphics[width=.45\hsize]{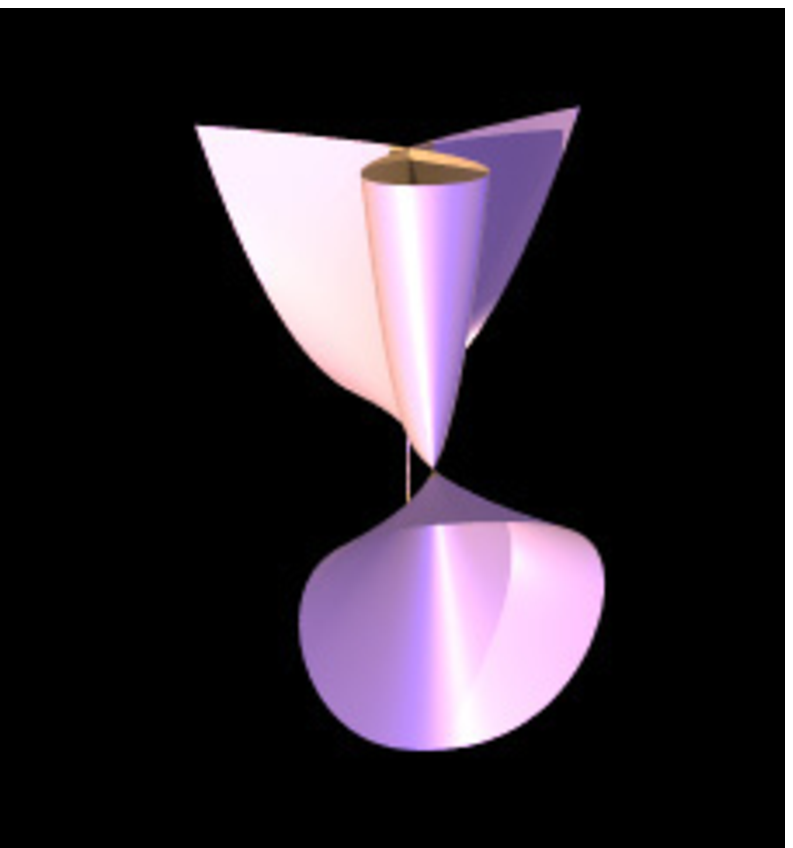}}}
\def\quasteB{{\includegraphics[width=.45\hsize]{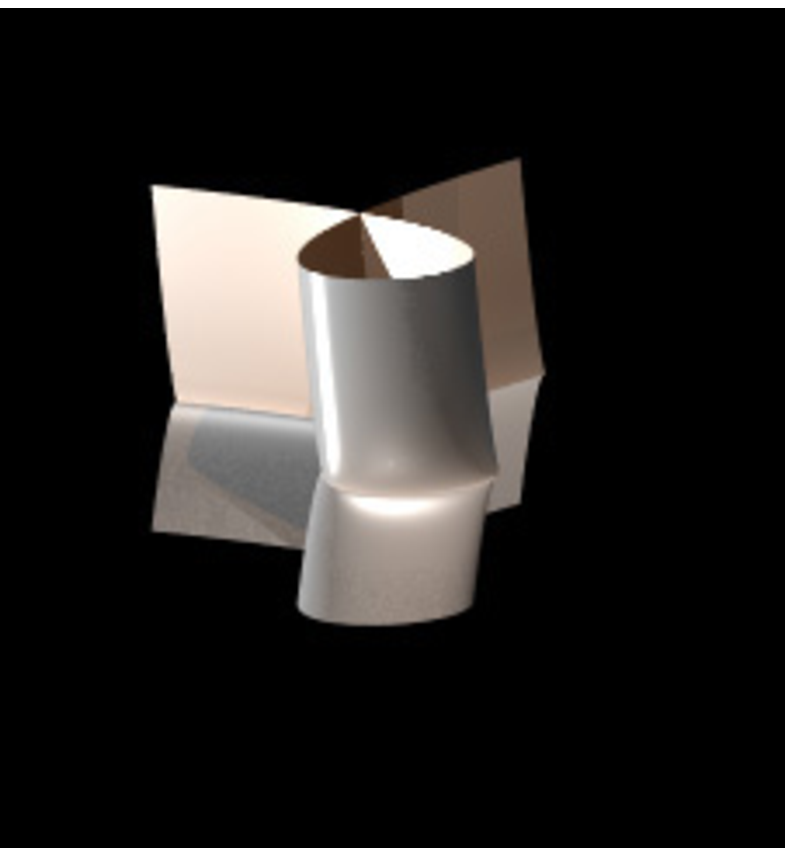}}}
\def\schneeflockeB{{\includegraphics[width=.45\hsize]{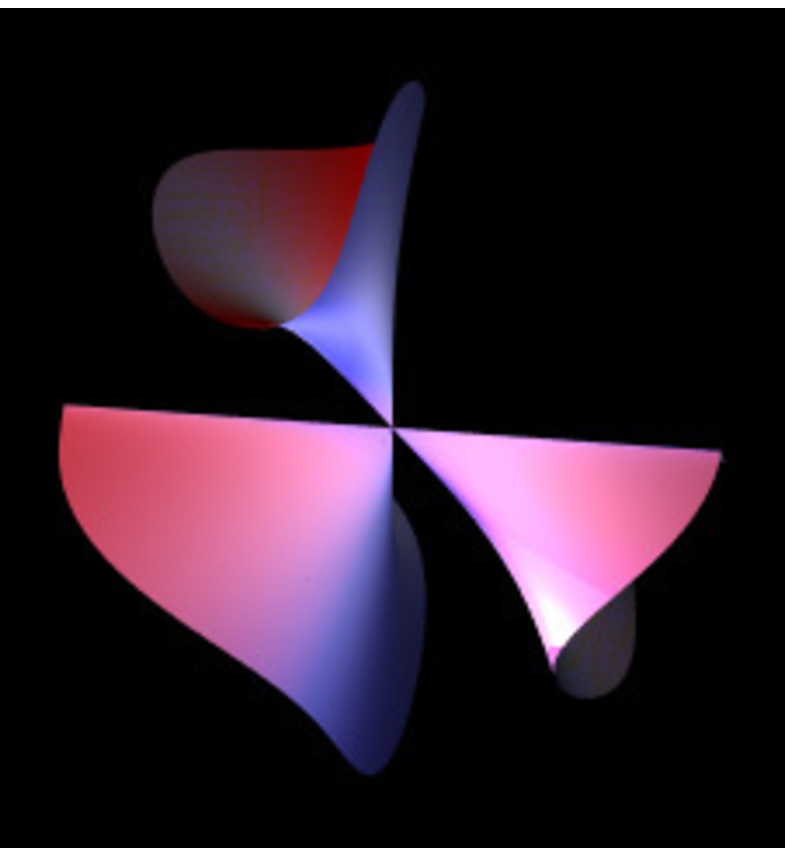}}}
\def\suessB{{\includegraphics[width=.45\hsize]{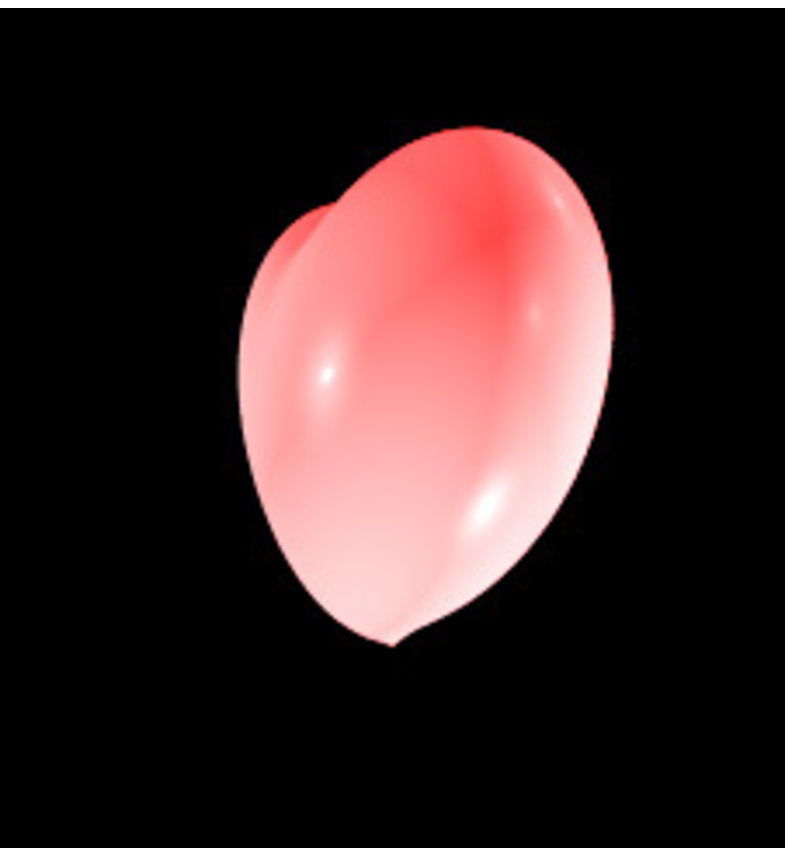}}}
\def\tanzB{{\includegraphics[width=.45\hsize]{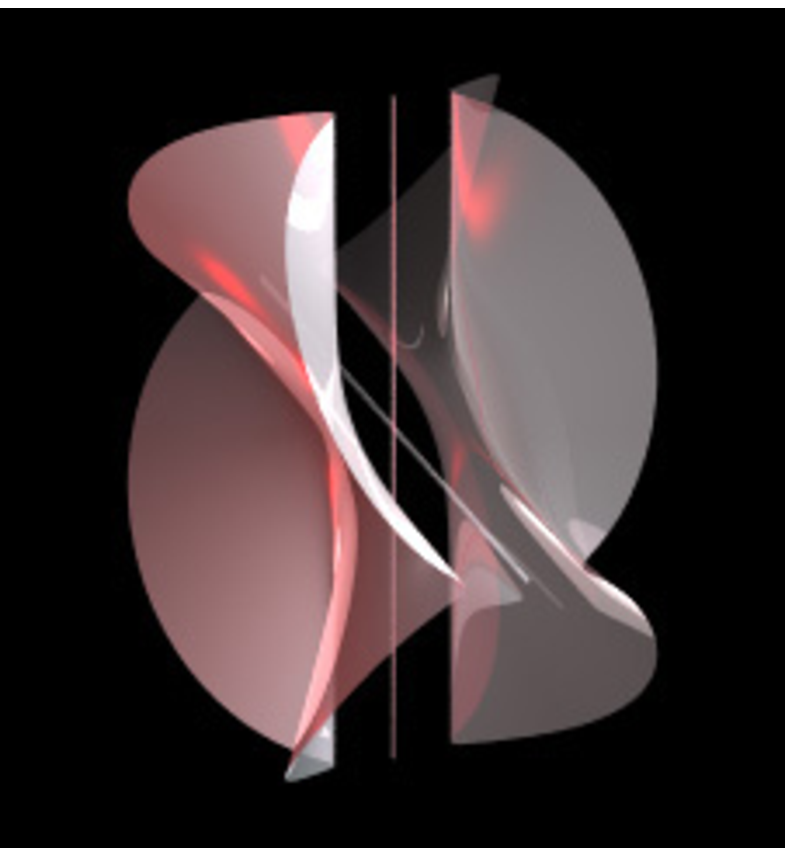}}}
\def\taubeB{{\includegraphics[width=.45\hsize]{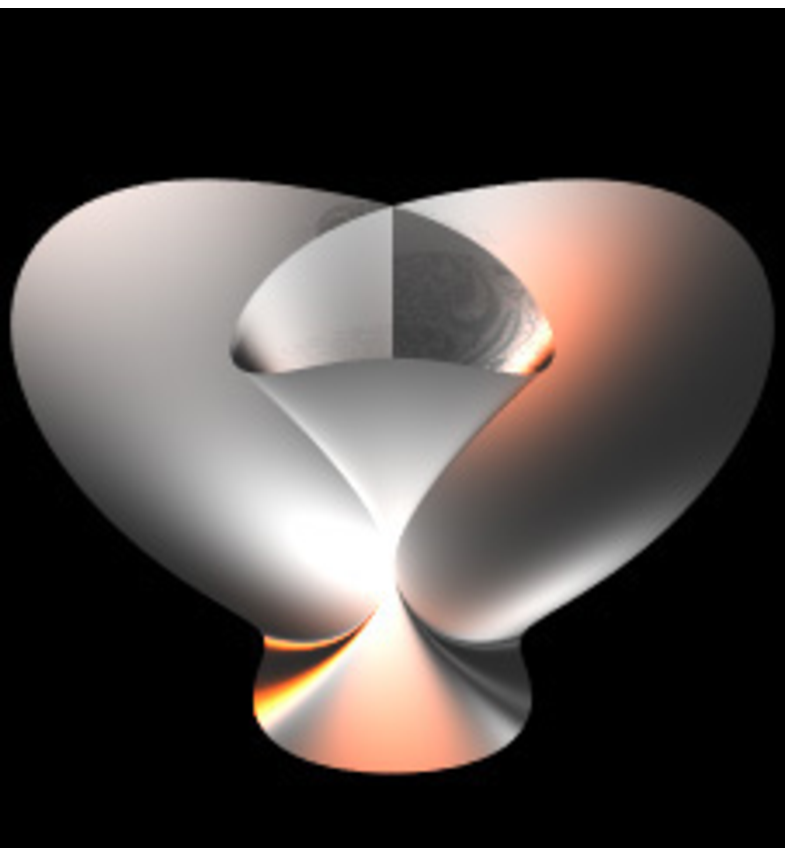}}}
\def\tuelleB{{\includegraphics[width=.45\hsize]{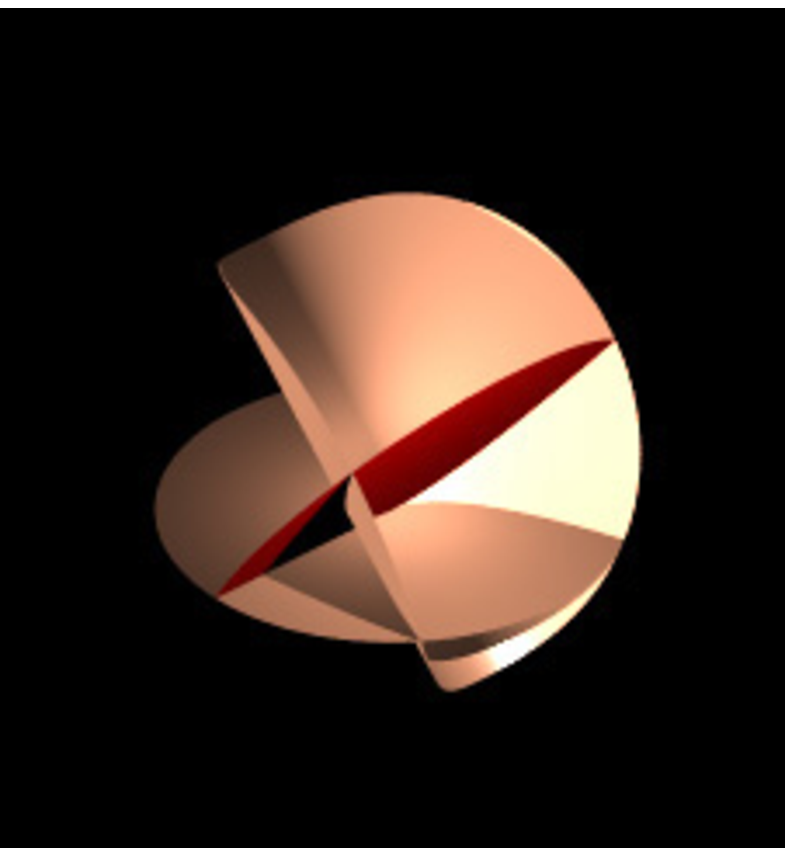}}}

\def\visavisB{{\includegraphics[width=.45\hsize]{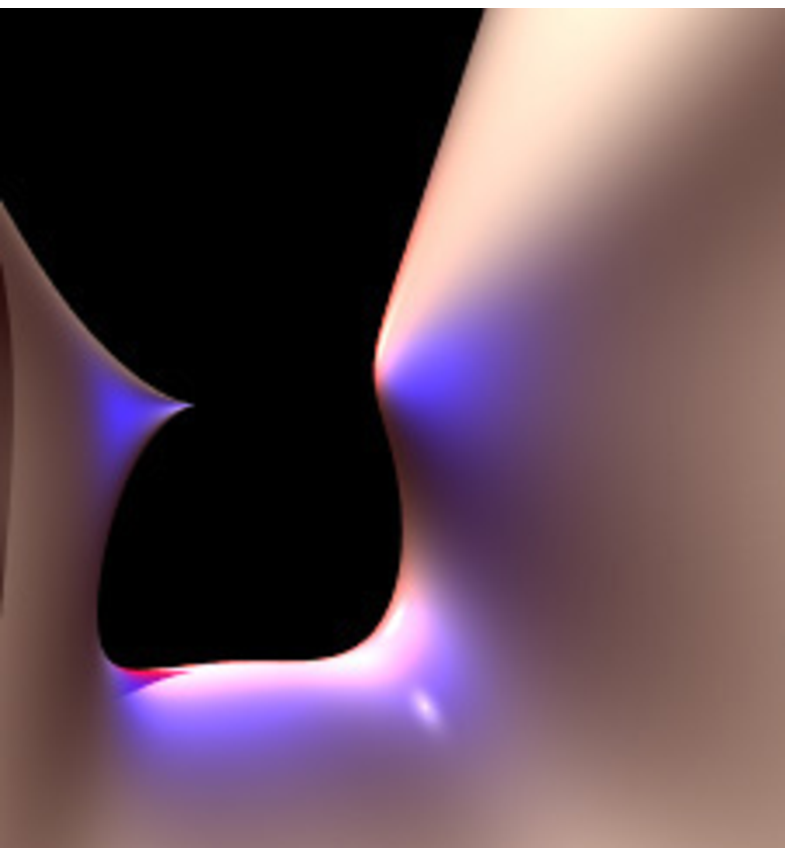}}}
\def\zeppelinB{{\includegraphics[width=.45\hsize]{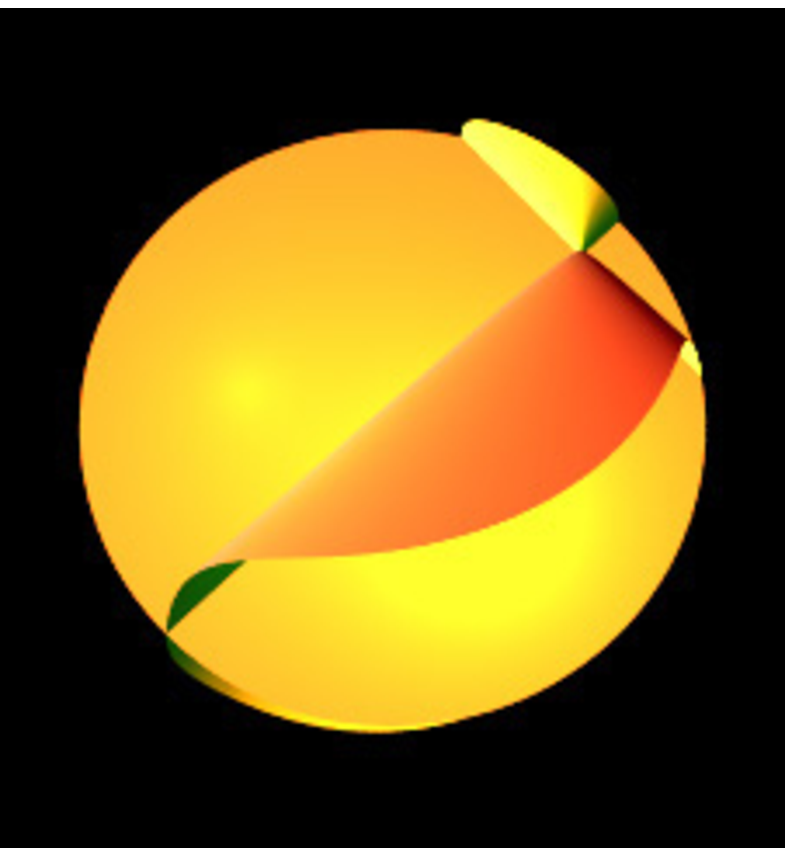}}}
\def\zitrusB{{\includegraphics[width=.45\hsize]{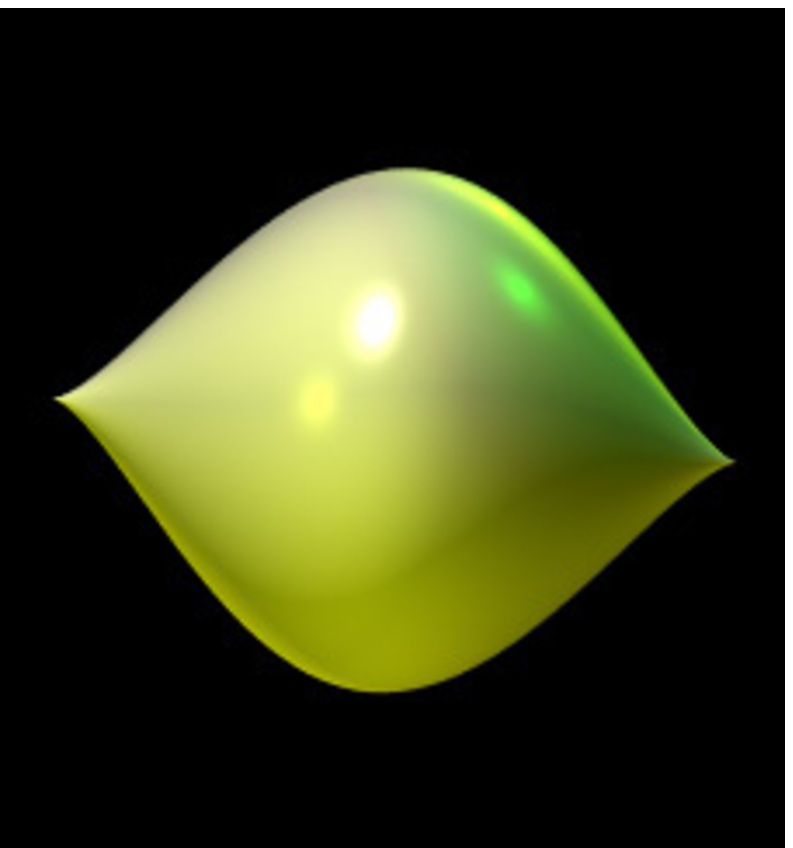}}}

\font\Times=ptmr at 12pt
\font\Timesten=ptmr at 10pt
\font\bf=ptmb at 12pt
\font\Bf=ptmb at 16pt

\font\it=ptmri at 12pt

\Times


\def\ltextindent#1{\hbox to \hangindent{#1\hss}\ignorespaces}
\long\def\ignore#1\recognize{}
\long\def\<#1\>{\vbox{\leftskip=1truecm{\eightrm #1}}}
\def\big{\bigskip}
\def\med{\medskip}

\def\hs{\hskip}
\def\vs{\vskip}
\def\cl{\centerline}

\def\{{\lbrace}
\def\}{\rbrace}

\def\R{{\Bbb R}}

\def\CC{{\cal C}}


 \cl{\Bf UFOs -- UNIDENTIFIED FIGURATIVE OBJECTS }\big

 \cl {\Bf A GEOMETRIC CHALLENGE}\big\big

 \cl{ C. BRUSCHEK, S. GANN, H. HAUSER, D. WAGNER, D. ZEILLINGER}\vs .5cm

 \cl { Institut f\"ur Mathematik}\par
 \cl {Universit\"at Innsbruck, Austria} \vs 1.3cm

\footnote{}{\Timesten e-mail: csac8572/sebastian.gann/herwig.hauser/csac8698/dominik.zeillinger@uibk.ac.at. \par

Supported in part by FWF-project P 15551.}

The title, though possibly provocative,  points at an identification task in algebraic geometry that is far from being trivial or resolved: the {\it Geometric Recognition Problem}. It  consists in detecting from given geometric objects -- in this case varieties -- their algebraic definition. More precisely, algebraic descriptions of objects which have the same visual properties as the given ones (there may not be uniqueness).\med

The article presents views of 24 real algebraic surfaces of degree mostly less than 7. The challenge is to find suitable equations, respectively parametrizations, for all of them. Of course, the pictures only allow $\CC^1$-interpretations, since the human eye cannot recognize the failure of $\CC^2$-differentiability. To indicate geometrically higher order features one had to depict the (iterated) Nash-modifications or transforms of the varieties under blowups.\med

Each picture comes with a short description of its main geometric configurations which may help to reveal the genesis. Already here we have to be somewhat vague -- there seems to be no generally accepted formal language to describe geometry. We do not propose here any such language, nor a procedure how to find algebraic equations for the surfaces. This is a vast field and could be the subject of further investigations. In this sense the article shall only serve as an appetizer: First, to contemplate geometric forms in order to find their exact mathematically rigorous description. Secondly, to search for algebraic procedures in order to construct surfaces from such a description using basic building blocks (e.g., by moving curves or by deforming toric surfaces).\med

We do not indicate appropriate equations and parametrizations for the surfaces which appear in this article. This is on purpose. Otherwise, the reader could quickly verify that the alge\-braic description is ok and pass to the next one. No! It is much more puzzling to try yourself finding an equation of an object which you believe to capture geometrically perfectly well.\med

The Geometric Recognition Problem may initiate new ideas and techniques to algebraic geometry, and not only over the reals or in dimension two. It may thus represent an (even though minor) counterbalance to the domination of algebraic, analytic and topological methods in geometry.\med

To satisfy the curiosity of the interested reader (not everybody will have time to determine the algebraic origin of the surfaces) we publish daily in december 2005 the equation or parametrization of {\it one} surface of the article in form of a calendar. The information can be found in the net at {\tt www.hh.hauser.cc}.
%
%
There, also a short animation of the surface (rotation, deformation, zooming, ...) is shown.\med

\goodbreak


Before passing to the series of pictures, we briefly describe the rules of the game. The coordinate system is in general fixed as follows: $x$ red, $y$ green, $z$ blue. \vs .7cm

\cl{\coordinates}\big\big

This may, however, vary according to the position of the camera.  The camera always looks at the origin -- in a few cases with slight deviations if required to have an authentic view. In the figures, the surfaces are clipped either with a ball or a box (mostly centered at the origin).\med

By {\it component} we shall always mean {\it visual component} (in a heuristic sense) with respect to real three-space and not in the algebraic sense as an irreducible component with respect to the Zariski topology. The same remark applies to the singular locus (the algebraic one and the one we see).
\med

As mentioned before, the pictures can only provide $\CC^1$-information. This is a drawback which produces a certain ambiguity when searching for equations. We therefore indicate occasionally additional properties which can help to focus the search on specific ranges of equations.\med

The pictures were produced by the authors with the ray-tracing program POV-Ray. One way to find the equations or parametrizations would be to scan the pictures, to choose a sufficiently narrow grid and to interpolate from the resulting reference points. This is not what we have in mind or what we are interested in. Our approach is synthetic: Find the algebraic description of the surface by understanding algebraically its geometric properties and construction rules and by then expressing them accurately in a suitable language. \med\med

{\it Achtung}: There are no results nor theorems in this paper, not even the description of a construction or an algorithm. This can be seen as an {\it affront} in comparison to the traditional way of publishing in mathematics. It could be! Nevertheless, the authors believe that mathematical research is based on observing and understanding interesting phenomena. The investigation is often more valuable (and satisfactory) than the result itself. Here, with this article, we wish to start such an investigation, not to conclude it.\med

It is amazing to see how mathematicians, when confronted with visualizations of surfaces they prove theorems about, are sometimes surprised to realize that the geometric object is indeed the subject of their investigation. We made this experience at several occasions. It shows that geometry and geometric contemplation have lost importance -- despite the fact that algebraic geometers try to describe and understand geometric phenomena.\med

The names of the surfaces were chosen by the authors. The reader will have no problems in finding the appropriate translation to her/his language.

 \vfill\eject

{\bf 1    \hs .2cm Helix:} degree 4 \med  

We start with a surface whose cross-sections with planes parallel (but different)  to the coordinate planes look like lemniscates (for $x=c$ and $z=c$) or unions of two hyperbolas (for $y=c$). The section with $y=0$ is the union of the $x$- and $z$-axis. The surface has reflections about the $xz$-plane and $90^\circ$ rotations around the $y$-axis as finite symmetries. More generally, one may ask for helices with an arbitrary number of petals.\big \big

\cl{\helixA \hs 1cm \helixB}\vs 1cm   


{\bf 2   \hs .2cm Tanz:} degree 4 \med  

In this surface, two two-dimensional (visual) components are opposed to each other by a reflection about the $xz$-plane. In addition, the $x$- and the $z$-axes form one-dimensional components. The intersections with the planes $x=\pm 1$ are translates of the $y$- and $z$-axis. The complexification of the surface coincides with the complexification of the Helix.\big\big

\cl{\tanzA \hs 1cm \tanzB}\vs 1cm\vfill\eject


{\bf 3   \hs .2cm Columpius:} degree 4\med  

Here, the original equation is difficult to find from inspection. In the picture, it is a perturbation of the equation $x^3y+y^3z+z^3x=0$ of the Klein Quartic producing three holes. The third hole is ``perpendicular'' to the two others and appears in the second figure. To relax the problem, one can ask for finding the equation of surfaces with a prescribed configuration of holes. \big \big

\cl{\columpiusA \hs 1cm \columpiusB}\vs 1cm


{\bf 4   \hs .2cm  Schneeflocke:}  degree 5 \med  

The surface is formed by  three components meeting at the origin. Two of them have the $y$-axis as singular locus and are generated by a contracting tacnode. The third component lies above the $z$-plane and has a funnel like shape. The surface contains the $y$- and $z$-axis and the diagonal $x=y+z=0$.\big\big

\cl{\schneeflockeA \hs 1cm \schneeflockeB}\vs 1cm\vfill\eject


{\bf 5   \hs .2cm Pirat:} degree 4 \med  

This is the surface one obtains when blowing up the the origin of the Whitney umbrella $x^2=y^2z$. At the meeting point of the one-dimensional component with the surface we have the local singularity of the Whitney umbrella (a vertical singular axis). The isolated singularity of the surface is an ordinary double point. \big\big

\cl{\piratA \hs 1cm \piratB}\vs 1cm


{\bf 6   \hs .2cm Vis-\`a-vis:} degree 4  \med  

Here, a cusp like isolated singularity is opposed to a local maximum (when considered in the direction of the horizontal axis). This surface is a perturbation of the surface Flirt further on. \big\big

\cl{\visavisA \hs 1cm \visavisB}\vs 1cm\vfill\eject


{\bf 7  \hs .2cm Dattel:} degree 2\med  

The equation of the ellipsoid is well known. \big\big

\cl{\dattelA \hs 1cm\dattelB}\vs 1cm


{\bf 8   \hs .2cm Nadel\"ohr:} degree 6 \med  

In contrast to the preceding, a suitable equation for the surface is hard to find. The main feature, aside the hole, is the self-tangency at the origin. Moreove, the right hand part has an isolated singular point there, created by the vanishing of a moving parabola.\big\big

\cl{\nadeloehrA \hs 1cm \nadeloehrB}\hs 1cm\vfill\eject


{\bf 9   \hs .2cm Cube:} degree 4 or 6 \med  

The equation of the rounded cube is relatively easy to determine, due to its symmetry properties. Taking an equation of degree 6 instead of 4 improves the approxiamtion with an exact cube.\big\big

\cl{\cubeA \hs 1cm \cubeB}\vs 1cm


{\bf 10  \hs .2cm T\"ulle:} degree 4 \med  

Two planes and one smooth surface intersect pairwise along the $x$-axis, respectively two parabolas, all three tangent to each other. This is an example of a non-normal crossings singularity with three components, where all pairwise intersections are smooth, but the common intersection is (scheme-theoretically) singular. In the second figure, the vertical plane has been tilted.\big\big

\cl{\tuelleA \hs 1cm\tuelleB}\vs 1cm\vfill\eject


{\bf 11   \hs .2cm Durchblick:} degree 4\med  

Another deformation of the Klein Quartic, but with just one hole (cf. with Columpius). The ondulation of the intersection of the surface with the sphere is typical. The cross-section of the hole is close to a tacnode (closed curve with one cusp singularity).\big\big

\cl{\durchblickA \hs 1cm \durchblickB}\vs 1cm 


{\bf 12   \hs .2cm  Lilie:} degree 5\med  

This Monet-like surface is hard to describe with words. The singularity at $0$ seems to be isolated (in the real picture), but the algebraic singular locus is more complicated.\big\big

\cl{\lilieA \hs 1cm \lilieB}\vs 1cm\vfill\eject


{\bf 13   \hs .2cm S\"uss:} degree 6 \med  

The equation of the heart can easily be found in the net.\big\big

\cl{\suessA\hs 1cm\suessB}\vs 1cm 


{\bf 14 \hs .2cm Quaste:} degree 9\med  

The cartesian product of the plane cusp with the node gives a surface in $\R^4$. Here, we see an embedding into $3$-space by taking a suitable projection. The singular locus is the union of a cusp and a node, meeting transversally at the origin. The parametrization is much simpler than the equation (the latter involves 24 monomials).\big\big

\cl{\quasteA \hs 1cm \quasteB}\vs 1cm \vfill\eject


{\bf 15   \hs .2cm Croissant:}   \med  

The pinched torus has a complicated equation, so that it is easier to determine a parametrization. The isolated singularity is cone-like (i.e., an ordinary double point). This surface is a typical example for illustrating characteristic classes and radial vector fields.\big\big

\cl{\croissantA \hs 1cm\croissantB}\vs 1cm 


{\bf 16   \hs .2cm Dromo:} degree 4  \med  

The symmetry of this surface is given by reflection about the $yz$-plane. The singularity at $0$ is similar to the one of the lilie (though simpler). The cross-section with the plane $y=0$ gives a plane curve which is a variation of the ordinary cusp $x^2=z^3$.\big\big

\cl{\dromoA\hs 1cm \dromoB}\vs 1cm\vfill\eject


{\bf 17   \hs .2cm Flirt:} degree 4 \med  

This surface is a deformation of Vis-\`a-vis. Two cusp like singularities meet with the same tangent at $0$. Globally, they form one hole, whose cross-section with the plane $x=0$ is (essentially) a circle.\big\big

\cl{\flirtA \hs 1cm \flirtB}\vs 1cm 


{\bf 18   \hs .2cm    Zeppelin:} degree 5 \med  

Along the $y$-axis, we have outside the origin a normal crossings singularity given by two transversally intersecting planes. At the origin, one of the components winds around (with tangent plane equal to the $xz$-plane), the other continues straight ahead.\big\big

\cl{\zeppelinA \hs 1cm \zeppelinB}\vs 1cm\vfill\eject


{\bf 19   \hs .2cm Daisy:} degree 6 \med  

The singular locus consists of two components, two cusps, meeting transversally at the origin. The two components of the surface have cone like shape taken over a closed curve with two cusp singularities (the intersection of the surface with a sphere). \big \big

\cl{\daisyA \hs 1cm \daisyB}\vs 1cm 


{\bf 20   \hs .2cm  Clip:} degree 5 \med  

Despite its modest appearance, the surface has an interesting algebraic background.
Its (algebraic) singular locus is the twisted space curve with parametrization $(t^3, t^4,t^5)$. The equation of our example has minimal degree with this property. The singular curve cannot be seen, because outside the origin  the surface is $\CC^1$-smooth along its (algebraic) singularities. \big\big

\cl{\clipA \hs 1cm \clipB}\vs 1cm \vfill\eject


{\bf 21  \hs .2cm Zitrus:} degree 6 \med  

Surface of revolution of a plane curve with two cusp like singularities obtained by rotation about the $y$-axis. There is a reflection symmetry with respect to the $xz$-plane.\big\big

\cl{\zitrusA \hs 1cm \zitrusB}\vs 1cm 

\big

{\bf 22 \hs .2cm Leopold:} degree 6 \med  

The ellipsoid can be deformed by squeezing it towards the origin outside the coordinate planes. This already gives a good hint how to find the equation of the present surface. It is symmetric with respect to the permutation of $x$, $y$ and $z$ (up to a homothety).\big\big

\cl{\leopoldA \hs 1cm \leopoldB}\vs 1cm \vfill\eject

\big


{\bf 23   \hs .2cm Taube:} degree 5 \med  

Famous example of singularity theory, the discriminant of the miniversal deformation of an $A_3$-singularity. The value of the coefficients play a decisive role. The second figure is obtained by changing one coefficient of the defining equation.\big\big

\cl{\taubeA \hs 1cm \taubeB}\vs 1cm 

\big

{\bf 24   \hs .2cm Distel:} degree 6 \med  

Similar construction as for Leopold, stretching a sphere towards the origin except at its six intersection points with the coordinate axes. In the figure, the perturbation coefficient was chosen rather big so as to produce the acute points.\big\big

\cl{\distelA \hs 1cm \distelB}\vs 1cm \vfill\eject

\big
\vfill\eject\end

\hrule \big


{\bf ++}   \hs .2cm Klein Q:\med  

\big


{\bf ++}   \hs .2cm  Michelangelo:\med  

\big

{\bf ++}   \hs .2cm T\"uchlein:\med  
\big

{\bf ++}   \hs .2cm E-7:\med  

\big

{\bf ++}   \hs .2cm Kegel mit Falte:\med  

\big


{\bf ++}   \hs .2cm Flapjack:\med  

\big


{\bf ++}   \hs .2cm I-c-e :\med  

\big


{\bf ++}    \hs .2cm Windkanal:\med  

\big


{\bf ++}    \hs .2cm Zuckerl:\med  

\big


{\bf ++}   \hs .2cm Unhold:\med  

\big


\vfill\eject \end